\RequirePackage[l2tabu, orthodox]{nag}
\documentclass[11pt]{amsart} %

\usepackage{amssymb,amscd,amsthm,mathrsfs}
\usepackage[utf8]{inputenc}
\usepackage[centertags]{amsmath}
\usepackage{graphicx}
\usepackage{setspace}
\usepackage{url}

\usepackage[capitalize]{cleveref}

\newcommand{\T}{\mbox{$\mathbb{T}$}}

\newcommand{\cG}{\mathcal{G}}

   \def\TT{{\mathbb T}}
   
 \def\ZZ{{\mathbb Z}}

\newtheorem*{thm*}{Theorem}
\newtheorem*{prop*}{Proposition}
\newtheorem*{cor*}{Corollary}
\newtheorem*{goal*}{Goal}

\newtheorem{thm}{Theorem}[section]

\newtheorem{quest}[thm]{Question} %
\newtheorem{cor}[thm]{Corollary}

\newtheorem{prop}[thm]{Proposition}

\newcommand{\bi}{\begin{itemize}}
\newcommand{\ei}{\end{itemize}}

\theoremstyle{definition}

\theoremstyle{remark}
\newtheorem{obs}[thm]{Remark}

\newtheorem{lem}[thm]{\bf Lemma}
\crefname{thm}{Theorem}{Theorems}
\Crefname{thm}{Theorem}{Theorems}
\crefname{figure}{Figure}{Figures}

\newcommand{\bbR}{\mathbb{R}}

\newcommand{\bbZ}{\mathbb{Z}}

\newcommand{\bbT}{\mathbb{T}}
\newcommand{\bbS}{\mathbb{S}}
\newcommand{\cE}{\mathcal{E}}

\newcommand{\Heis}{\mathcal{H}}

\newcommand{\cN}{\mathcal{N}}

\newcommand{\Es}{E^s}
\newcommand{\Ec}{E^c}
\newcommand{\Eu}{E^u}
\newcommand{\Ecu}{E^{cu}}
\newcommand{\Ecs}{E^{cs}}

\newcommand{\Ecsf}{\Ecs_f}

\newcommand{\EsA}{\Es_A}
\newcommand{\EcA}{\Ec_A}
\newcommand{\EuA}{\Eu_A}

\newcommand{\EcsA}{\Ecs_A}

\newcommand{\Esg}{\Es_g}

\newcommand{\Eug}{\Eu_g}

\newcommand{\cW}{\mathcal{W}}
\newcommand{\Ws}{\cW^s}

\newcommand{\Wu}{\cW^u}

\newcommand{\tW}{\tilde{\cW}}
\newcommand{\tWs}{\tW^s}
\newcommand{\tWc}{\tW^c}
\newcommand{\tWu}{\tW^u}

\newcommand{\sig}{\sigma}
\newcommand{\Sig}{\Sigma}

\newcommand{\loc}{\operatorname{loc}}
\newcommand{\graph}{\operatorname{graph}}

\newcommand{\cF}{\mathcal{F}}
\newcommand{\Fc}{\cF^c}
\newcommand{\Fcu}{\cF^{cu}}
\newcommand{\Fcs}{\cF^{cs}}

\newcommand{\Fcsep}{\cF^{cs}_{\epsilon}}

\newcommand{\cA}{\mathcal{A}}

\newcommand{\Acu}{\cA^{cu}}
\newcommand{\Acs}{\cA^{cs}}
\newcommand{\Aus}{\cA^{us}}

\newcommand{\tF}{\tilde{\cF}}

\newcommand{\tFcu}{\tF^{cu}}
\newcommand{\tFcs}{\tF^{cs}}

\newcommand{\tFcsep}{\tF^{cs}_{\epsilon}}

\newcommand{\tA}{\tilde{\cA}}
\newcommand{\tAs}{\tA^s}
\newcommand{\tAc}{\tA^c}
\newcommand{\tAu}{\tA^u}
\newcommand{\tAcu}{\tA^{cu}}
\newcommand{\tAcs}{\tA^{cs}}
\newcommand{\tAus}{\tA^{us}}

\newcommand{\Lcs}{L^{cs}}
\newcommand{\Lcu}{L^{cu}}

\newcommand{\cL}{\mathcal{L}}

\newcommand{\tL}{\tilde{L}}
\newcommand{\tM}{\tilde{M}}
\newcommand{\tf}{\tilde{f}}
\newcommand{\tg}{\tilde{g}}
\newcommand{\tilh}{\tilde{h}} %

\newcommand{\tvp}{\tilde{\varphi}}

\newcommand{\ep}{\epsilon}
\newcommand{\lam}{\lambda}
\newcommand{\Lam}{\Lambda}
\newcommand{\gam}{\gamma}
\newcommand{\Gam}{\Gamma}
\newcommand{\tgam}{\tilde{\gam}}

\newcommand{\Diff}{\operatorname{Diff}}

\newcommand{\dist}{\operatorname{dist}}
\newcommand{\length}{\operatorname{length}}
\newcommand{\volume}{\operatorname{volume}}

\newcommand{\supp}{\operatorname{supp}}

\newcommand{\Hcs}{H^{cs}}
\newcommand{\Hcu}{H^{cu}}

\newcommand{\Mane}{Ma\~n\'e}
\newcommand{\Diaz}{D\'iaz}

\newcommand{\qandq}{\quad \text{and} \quad}

\newcommand{\inv}{^{-1}}

\newcommand{\TLM}{T_\Lam M}

\newenvironment{idea}[1][\unskip]
{\begin{proof}[Idea of proof #1]}
{\end{proof}}

\theoremstyle{remark}

\providecommand{\acknowledgement}{{\noindent\bf Acknowledgements}\quad}

\newcommand{\en}{\subset}
\newcommand{\trans}{\mbox{$\,{ \top} \;\!\!\!\!\!\!\raisebox{-.3ex}{$\cap$}\,$}}

\author[A. Hammerlindl]{Andy Hammerlindl}
\address{School of Mathematical Sciences,
     Monash University,
      Victoria 3800 Australia}
\urladdr{\url{http://users.monash.edu.au/~ahammerl/}}
\email{andy.hammerlindl@monash.edu}

\author[R. Potrie]{Rafael Potrie}
\address{CMAT, Facultad de Ciencias, Universidad de la Rep\'ublica, Uruguay}
\urladdr{\url{www.cmat.edu.uy/~rpotrie}}\email{rpotrie@cmat.edu.uy}

\title[Partial hyperbolicity and classification]{Partial hyperbolicity and classification: a survey.}
\thanks{A.H. was partially supported by ARC Grant DP120104514.
R.P. was partially supported by CSIC group 618.}

\begin{document}
\begin{abstract}
  This paper surveys recent results on classifying partially hyperbolic
  diffeomorphisms.
  This includes the construction of branching foliations
  and leaf conjugacies on three-dimensional manifolds with solvable
  fundamental group.
  Classification results in higher-dimensional settings
  are also discussed.
  The paper concludes with an overview of the construction of new partially
  hyperbolic examples derived from Anosov flows. 
  \bigskip


\end{abstract}

\maketitle

\tableofcontents

\section{Introduction} %

This paper surveys the work in recent years on the classification problem for
partially hyperbolic dynamical systems.
One motivation for research into partial hyperbolicity is the study of
invariant objects in dynamical systems which persist under perturbation.
To this end, we start the survey by defining isolated sets and investigating
their properties under changes to the dynamics.

\subsection{Persistent manifolds}
Consider a diffeomorphism $f:M \to M$ on a Riemannian manifold.
A compact set $\Lam$ is 
\emph{invariant}
if $f(\Lam) = \Lam$.
It is
\emph{isolated}
if there is a neighborhood $U \supset \Lam$ such that
$\Lam = \bigcap_{n \in \bbZ} f^n(U).$
Note that any isolated set is automatically invariant.
Let $\Diff^r(M)$ denote the space of diffeomorphisms $f:M \to M$   
equipped with the $C^r$ topology \cite{Hirsch}.
For an isolated set $\Lam = \Lam_f$
and $g \in \Diff^r(M)$ close to $f$, the
\emph{continuation} of $\Lam_f$ is defined by
$\Lam_g = \bigcap_{n \in \bbZ} g^n(U)$.

In general, the continuation $\Lam_g$ could be very different in structure
from the original $\Lam_f$.
The simplest example is where $\Lam_f$ is a singleton set $\Lam_f = \{x\}$.
One can verify that $\Lam_g$ is still a singleton set for all nearby $g$
exactly when the derivative $Df_x:T_x M \to T_x M$ is hyperbolic
(i.e., with no eigenvalues of modulus one).

Consider now the case where the invariant set $\Lam$ is a submanifold.
Let $T_\Lam M := \{T_x M: x \in \Lam\}$ denote the tangent bundle $TM$
restricted to $\Lam$.
A \emph{subbundle} $E \subset T_\Lam M$ is a continuous choice of subspace
$E(x) \subset T_x M$ for each $x \in \Lam$.
The subbundle is \emph{invariant} if $Df(E) = E$ where $Df:TM \to TM$
is the derivative of $f$.
For example, if $\Lam  = f(\Lam) \subset M$ is an embedded submanifold,
then $T \Lam$, the set of vectors tangent to $\Lam$,
is an invariant subbundle of $T_\Lam M$.
An invariant submanifold $\Lam$ is
\emph{normally hyperbolic}
if there is a splitting of $T_\Lam M$ into invariant subbundles
\[
    T_\Lam M = \Es \oplus T \Lam \oplus \Eu
\]
and $n  \ge  1$
such that at each $x \in \Lam$
the inequalities
\[
    \|Df^n v^s\| < 1 < \|Df^n v^u\|
    \qandq
    \|Df^n v^s\| < \|Df^n v^c\| < \|Df^n v^u\|
\]
hold
for all unit vectors $v^s \in \Es(x)$, $v^c \in T_x \Lam$ and $v^u \in \Eu(x)$.
In words,
this means that the expansion and contraction normal to the submanifold
dominates any expansion or contraction tangent to the submanifold.

An isolated $C^1$-submanifold $\Lam_f$ is called
\emph{persistent}
if $\Lam_g$ is also a $C^1$-submanifold close to the original $\Lam_f$
for all $g$ in a $C^1$ neighborhood of $f$.

\begin{thm} \label{thm:nhper}
    Normally hyperbolic submanifolds are persistent.
\end{thm}

This result, in one form or another, has been proved many times over,
with the techniques of the proof
going back at least to a 1901 paper of Hadamard on stable manifolds.
See \cite[Section 1]{HPS} for references.

In his 1973 thesis, R.~\Mane{} \cite{ManheThesis} proved the converse.

\begin{thm}
    [\Mane{}] \label{thm:pernh}
    Persistent submanifolds are normally hyperbolic.
\end{thm}
Analogous results also hold for continuous-time dynamical systems
 \cite{Osipenko}. %

\medskip{}

In their 1977 book, Hirsch, Pugh, and Shub
developed extensions to \cref{thm:nhper} in great detail \cite{HPS}.
One case considered is when $\Lam$ supports a lamination.
That is, $\Lam$ decomposes as a disjoint union of submanifolds,
called leaves,
such that if $L_x$ denotes the leaf through $x \in \Lam$,
then the tangent space $T_x L_x \subset T_x M$
depends continuously on $x$.
This lamination is invariant if both $f(\Lam) = \Lam$
and $f(L_x) = L_{f(x)}$ for all $x \in \Lam$.
We use $\cF$ to denote the lamination and say that $\Lam_f$ is the
support of $\cF$.
Normal hyperbolicity is then defined as above with $T_x L_x$ in place
of $T_x \Lam$.

Hirsch, Pugh, and Shub proved a limited form of persistence
for normally hyperbolic laminations \cite[Theorem 6.8]{HPS},
but a priori the lamination structure may be destroyed under perturbation.
For instance, if $L_x$ and $L_y$ are distinct leaves in $\cF$,
then the corresponding submanifolds in $\Lam_g$ may a priori intersect
without coinciding.
If we add one additional assumption called ``plaque expansiveness,''
then these thorny issues disappear.

As it is technical, we omit the definition of plaque expansiveness.
To state the theorem, however, we define one more notion.
Suppose $f,g:M \to M$ are diffeomorphisms with
respective invariant laminations $\cF$ and $\cG$.
A \emph{leaf conjugacy}
between $(f, \cF)$ and $(g, \cG)$ is a homeomorphism
$h: \supp \cF \to \supp \cG$
such that if $L_x$ is a leaf of $\cF$, then $h(L_x)$
is a leaf of $\cG$ and
$h f (L_x) = g h (L_x)$.
That is, conjugacy holds on the level of leaves.

\begin{thm} \label{thm:peconj}
    Suppose $f:M \to M$ is a diffeomorphism and
    $\cF$ is an $f$-invariant, normally hyperbolic, plaque expansive
    lamination.
    Then, there is a $C^1$-neighborhood $\mathcal{U}$ of $f$
    such that for every $g \in \mathcal{U}$ there is a
    $g$-invariant, normally hyperbolic, plaque expansive lamination $\cG$
    and $(g, \cG)$ is leaf conjugate to $(f, \cF)$.
\end{thm}
Plaque expansive holds in every known example.

\begin{quest}
    Is every normally hyperbolic lamination also plaque expansive?
\end{quest}
\begin{quest} 
    What are all possible normally hyperbolic laminations?
    Can they be classified up to some notion of equivalence?
\end{quest}

In light of \cref{thm:peconj} the most natural notion of equivalence would be
leaf conjugacy.  This problem may be viewed as the major motivation
of the research detailed in this survey paper.
As the problem in its full generality is way too difficult at present,
we look at it in specific settings.

\subsection{Uniform hyperbolicity}

Normal hyperbolicity builds on the notion of uniform hyperbolicity
developed by Anosov, Smale, and others in the 1960s.
Consider a diffeomorphism $f:M \to M$.
A compact invariant subset $\Lam$ is
\emph{hyperbolic}
if there is an invariant splitting
\[
    T_\Lam M = \Es \oplus \Eu
\]
and $n  \ge  1$
such that
\[
    \|Df^n v^s\| < 1 < \|Df^n v^u\|
\]
for all unit vectors $v^s \in \Es$ and $v^u \in \Eu$.
A hyperbolic set may be viewed as a special case of a normally hyperbolic
lamination,
where the ``leaf'' through a point $x$ is just the singleton set $\{x\}$.
The leaf conjugacy in \cref{thm:peconj} is then a true conjugacy 
between the homeomorphisms $f|_{\Lam_f}$ and $g|_{\Lam_g}$.
Further, ``plaque expansiveness'' is just ``expansiveness'' here,
and is always known to hold for hyperbolic sets.

When the manifold $M$ is a hyperbolic set, the diffeomorphism
$f:M \to M$ is called \emph{Anosov}.
Several important classification results have been established for
Anosov diffeomorphisms and we discuss these in \cref{sec:classanosov}.

\subsection{Partial hyperbolicity}
Shortly after the development of normal
hyperbolicity, Brin and Pesin
developed a related notion of partial hyperbolicity \cite{BrinPesin}.
There are several ``flavours'' of partial hyperbolicity.
As this is a survey on the subject, we give a list of the various definitions.

Consider again a diffeomorphism $f:M \to M$ and assume there is a compact 
invariant set $\Lam$
with invariant splitting
\[
    \TLM = \Es \oplus \Ec \oplus \Eu
\]
such that at each $x \in \Lam$
the inequalities
\begin{equation} \label{equation:phineq}
    \|Df^n v^s\| < 1 < \|Df^n v^u\|
    \qandq
    \|Df^n v^s\| < \|Df^n v^c\| < \|Df^n v^u\|
\end{equation}
hold
for all unit vectors $v^s \in \Es(x)$, $v^c \in \Ec(x)$ and $v^u \in \Eu(x)$.
We assume $\Ec$ is non-zero, as otherwise
the set would be uniformly hyperbolic.
If exactly one of the bundles $\Es$ and $\Eu$ is non-zero,
then $\Lam$ is a
\emph{weakly partially hyperbolic set}.
If $\Es$, $\Ec$, and $\Eu$ are all non-zero,
then $\Lam$ is a
\emph{strongly partially hyperbolic set}.
The bundles $\Es$, $\Ec$, and $\Eu$ are called the
``stable'', ``center'', and ``unstable'' bundles respectively.

Note that any strongly partially hyperbolic set may also be viewed as
a weakly partially hyperbolic by, say, grouping the stable and center
bundles together into a single center bundle.

The definition of partial hyperbolicity above is pointwise;
the inequalities \eqref{equation:phineq}
hold independently at each $x \in \Lam$.
A partially hyperbolic set is \emph{absolutely} partially hyperbolic
if the inequalities \eqref{equation:phineq} hold for any points
$x,y,z \in \Lam$
and unit vectors $v^s \in \Es(x)$, $v^c \in \Ec(y)$, and $v^u \in \Eu(z)$.

Recently, Rodriguez Hertz, Rodriguez Hertz, and Ures constructed interesting
examples of pointwise partially hyperbolic systems that are impossible in the
absolute setting \cite{HHU}.

Throughout this paper, ``partially hyperbolic'' means 
strongly pointwise partially hyperbolic
unless otherwise noted.
We also primarily consider dynamical systems where the partially
hyperbolic set is the entire manifold $M = \Lam$ and say that the
diffeomorphism $f:M \to M$ is partially hyperbolic.

In this setting,
one can ask if there is an invariant foliation $\Fc$ tangent to the center
direction $\Ec$.  Such a foliation would necessary be normally hyperbolic.
We call a partially hyperbolic diffeomorphism
\emph{dynamically coherent}
if there are invariant foliations
$\Fcs$, $\Fcu$, and $\Fc$
tangent respectively to
$\Ecs$, $\Ecu$, and $\Ec$.

\begin{quest}
    Which partially hyperbolic systems are dynamically coherent?
\end{quest}
\begin{quest}
    Can the dynamically coherent systems be classified up to leaf conjugacy?
\end{quest}
\begin{quest}
    Can the non-dynamically coherent systems be classified up to some
    other notion of equivalence?
\end{quest}

\subsection{Robust transitivity}
Partially hyperbolic dynamics arise naturally when studying robustly transitive systems. 

An invariant set is
\emph{transitive}
if it contains a dense orbit.
That is, if there is $x \in \Lam$ such that
\[
    \Lam = \overline{ \{f^n(x) : n \in \bbZ \} }.
\]
An isolated set is
$C^r$-\emph{robustly transitive}
if the continuation $\Lam_g$ is transitive for all $g$
in a $C^r$-neighborhood of $f$.
Transitive hyperbolic sets are robustly transitive
and form a large class of examples.
\Mane{} showed that in dimension 2, the two notions are equivalent \cite{ManheContributions}.

\begin{thm}
    Suppose $M$ is a surface and $f:M \to M$ a diffeomorphism.
    A transitive set $\Lam$ is $C^1$-robustly transitive if and only if
    it is hyperbolic.
\end{thm}

Higher dimensions are more complicated.
In 1968, Shub constructed a robustly transitive example
on $M = \Lam = \bbT^4$
which is partially hyperbolic but not hyperbolic (see \cite[Section 8]{HPS}).
In 1978, \Mane{} constructed a similar example on $\bbT^3$ \cite{ManheContributions}.
In 1995, Bonatti and \Diaz{} developed the idea of a ``blender''
to construct many examples of invariant sets that are both
partially hyperbolic and robustly transitive \cite{BD-blender}.
In 1999, \Diaz{}, Pujals, and Ures proved the following theorem
in the converse direction \cite{DPU}.

\begin{thm}
    For a diffeomorphism on a 3-dimensional manifold
    $M$, any $C^1$-robustly transitive set must be
    weakly partially hyperbolic.
\end{thm}
\begin{quest}
    What are the robustly transitive sets in dimension 3?
    In particular, if $M$ is a compact 3-manifold $M$,
    then for which diffeomorphisms
    $f:M \to M$ is the entire manifold $M$
    a robustly transitive set?
    Can any 3-manifold $M$ support such a diffeomorphism $f$?
\end{quest} 

Robustly transitive sets in dimensions 4 and higher
must also have some form of invariant splitting \cite{BDP}
called ``volume partial hyperbolicity'' or just ``volume hyperbolicity.''
This notion is similar in form to partial hyperbolicity,
but there are robustly transitive sets in these higher dimensions
which are not strictly partially hyperbolic. See \cite[Chapter 7]{BDV}. 

\medskip{}

Interesting results also hold in the continuous time setting.
For flows in dimension 2,
the only robustly transitive sets are equilibrium points and attracting and
repelling orbits \cite{Peixoto}.

For flows in dimension 3, every $C^1$-robustly transitive set must be
weakly partially hyperbolic.
Examples include the Lorenz attractor, and
any robustly transitive set with an equilibrium point can, in some sense,
be viewed as a finite number of geometric Lorenz attractors ``glued together.''
For further details and precise statements see \cite{AraujoPacifico} or \cite[Chapter 9]{BDV} and references therein. 
\subsection{Stable ergodicity}

Finally,
we mention stable ergodicity as a major motivation for studying
partially hyperbolic dynamics.
Invariant measures often arise naturally in dynamics, for instance,
when studying symplectic dynamical systems
and specifically in the case of Hamiltonian systems
due to Liouville's theorem.
Consider a diffeomorphism $f$ of a compact manifold $M$
and suppose there is an invariant probability measure $\mu$ on $M$.
This pair $(f,\mu)$ is
\emph{ergodic}
if
\[
    \lim_{N \to \infty} \frac{1}{N}\sum_{k = 1}^N \phi(f^k(x))
    = \int_M \phi \, d \mu
\]
for almost every $x$.
In other words, $\mu$ captures the average behaviour of
almost every orbit of $f$.
For this discussion, we assume $\mu$ is equivalent to Lebesgue measure on $M$.
The pair $(f,\mu)$ is
\emph{stably ergodic}
if $(g,\mu)$
is also ergodic for every $\mu$-preserving diffeomorphism
$g$ close to $f$.
(For reasons lost to history, the term ``stably ergodic''
has won out over ``robustly ergodic.'')

In the 1960s,
D.~V.~Anosov demonstrated that volume-preserving Anosov diffeomorphisms
are stably ergodic \cite{Anosov}.  
This is a large part of the reason that they are called
``Anosov'' diffeomorphisms.
For decades, Anosov systems were the only known examples
of stable ergodicity, though some systems were known to be stably ergodic among certain interesting classes of perturbations.
A notable case is that of frame flows in negative curvature which represent a
class of genuinely partially hyperbolic flows, which under certain conditions
can be shown to be stably ergodic among frame flows
\cite{BrinPesin,Brin-Argument,BrinGromov}. 

In 1995, Grayson, Pugh, and Shub gave an example of
a partially hyperbolic diffeomorphism
which is stably ergodic \cite{GPS}.%
This lead to a frenzy of research in the subsequent years
showing that many volume-preserving partially hyperbolic systems
are stably ergodic.
See \cite{BuPSW,CHHU,HHU-Survey,PughShub,WilkinsonSurvey} for surveys of the results relating partial hyperbolicity
and stable ergodicity.  
See also \cite{HamUres,HamErgodic} which directly apply
classification and leaf conjugacy results
to the study of ergodicity.

\section{Classification and statements of the surveyed results} %

Many surveys on partial hyperbolicity already exist
\cite{BDV,BuPSW,HHU-Survey,hasselblattpesin,PughShub,WilkinsonSurvey}.
Particularly, a recent survey treats partially hyperbolic diffeomorphisms in
dimension 3 including a considerable section devoted to classification and
integrability \cite{CHHU}.

It is therefore a big challenge to write a survey
which complements these existing materials.
We felt that there were still
interesting things to be said, in particular regarding the classification
problem.
In this paper, many relevant aspects are omitted
or only briefly discussed if they are covered in detail in other surveys.
For instance,
it is standard when surveying a mathematical subject to start with a
list of known examples.
Instead, we refer the reader to the above surveys and to
\cite[Section 3]{CrovisierPotrie}.

\medskip

A geometric structure
such as a splitting of the tangent bundle into subbundles will likely impose
some form of restriction both on the manifold
and on the isotopy class of a diffeomorphism which leaves it invariant.
For Anosov diffeomorphisms in dimensions 4 or higher, however, the
classification problem is yet completely open.
For Anosov flows, the same problem arises.
Even in dimension 3, where many results are known, there
are still many open questions.
For example,
infinitely many hyperbolic 3-manifolds are known to support
Anosov flows \cite{Fenley} and infinitely many do not \cite{RSS},
but there is still an infinite class of hyperbolic manifolds where the question is wide
open.

Classifying partially hyperbolic systems is, in principle, much more
complicated and includes the classification of Anosov diffeomorphisms\footnote{Notice that even if there are Anosov systems where one can not choose a center bundle, their classification will be needed when studying partially hyperbolic skew products.} and
flows as subproblems. Indeed, the understanding of which manifolds and
isotopy classes admit partially hyperbolic systems is much less advanced than
for Anosov systems.

The bulk of this survey concerns the classification of partially hyperbolic diffeomorphisms in closed 3-dimensional manifolds. Higher
dimensions are discussed only briefly at the end. The main point, suggested
informally by Pujals in a conference and later formalized by Bonatti and
Wilkinson \cite{BW} is that even if classification of Anosov systems may be
out of reach, it is reasonable to attempt to \emph{compare} partially
hyperbolic diffeomorphisms to Anosov systems.
Many conjectures and perspectives were triggered by this approach and we refer to the recent survey
\cite{CHHU}.

Here we present a series of recent results and examples related to this
approach.
Of course, the choice of results is biased by the interests of the authors with a special emphasis on presenting both our previous and
ongoing joint work. 

\smallskip

There have
been essentially two approaches to the classification problem.
The first, initiated in \cite{BW}, looked at the behavior of the center
foliation in known examples and showed that this behavior characterized
the examples.
This has been quite successful and it has even lead to
progress in higher dimensions
\cite{bohnet2015quotient,carrasco2011compact,
gogolev2011compact,Hammerlindl,potrietrapping}
and it is partially surveyed
in \cref{sec:others}.

A different approach, initiated by the work of Brin, Burago and Ivanov
\cite{BBI,BI}, instead makes assumptions on the underlying topology of the
manifold and, in cases, also on the isotopy class of the diffeomorphism.

\smallskip
In the rest of this section, we restrict to the 3-dimensional case.
The study of 3-manifolds is quite advanced, and it makes sense to separate the
study of partially hyperbolic diffeomorphisms based on the geometry of the
manifold.

First, we consider ``small'' manifolds where the volume of a ball in the
universal cover is polynomial in the radius of the ball.
Equivalently, the fundamental group has polynomial growth (see \cite{GromovNil}).
In this setting, unstable curves grow exponentially fast under the dynamics
and in some sense wrap around themselves in the compact manifold.
This idea of curves in manifold ``wrapping around themselves''
under iteration may be made rigorous by considering the action in homology.
This was formalized and proved in \cite{BI} and extended in \cite{Parwani} to show the following: 

\begin{thm} \label{thm:phpoly}
Let $f$ be a partially hyperbolic diffeomorphism supported on a 3-manifold
$M$ such that the fundamental group has polynomial growth.
Then, (up to finite cover)
$M$ is a circle bundle over a torus,
and the action $f_*$ on homology $H^1(M,\bbR)$
has eigenvalues $\{\lam_i\}$ with $\min |\lam_i| < 1 < \max |\lam_i|$.
\end{thm}

\Cref{sec:integrability,sec:reebless} give an outline of results
featuring in the proof of \cref{thm:phpoly}.

Circle bundles over the torus may also be viewed as quotients
of nilpotent Lie groups by a discrete lattice.
That is, they are nilmanifolds.
Further, for every nilmanifold in dimension three
there is an automorphism of the nilpotent Lie group
which quotients down to a partially hyperbolic
diffeomorphism on the nilmanifold \cite[Appendix A]{HP2}.

Using this, the authors were able to classify all partially hyperbolic
diffeomorphisms on these manifolds \cite{HP}.
(See also \cite{Hammerlindl,HNil,Pot} for previous results.)

\begin{thm} \label{thm:nil}
Let $f$ be a partially hyperbolic diffeomorphism supported on a 3-manifold
$M$ such that the fundamental group is nilpotent
and suppose there is no compact surface tangent to $\Ecs$ or $\Ecu$.
Then, (up to finite cover)
$f$ is leaf conjugate to an algebraic map on $M$.
\end{thm}

\smallskip

We next consider the case where the fundamental group has exponential growth.
This comprises a much larger class of 3-manifolds.
Further, exponential growth is a necessary condition
for a manifold to support an Anosov flow \cite{PlanteThurston}.
In this class of ``big'' manifolds, the panorama is less clear. 
In this setting, the simplest manifolds to study are those with solvable
fundamental group.
These are the manifolds produced by the suspension of an Anosov
diffeomorphism, and the authors showed that the following classification
result holds \cite{HP2}.

\begin{thm} \label{thm:sol}
Let $f$ be a partially hyperbolic diffeomorphism supported on a 3-manifold
$M$ such that the fundamental group is solvable and has exponential growth.
Further, suppose there is no compact surface tangent to
$\Ecs$ or $\Ecu$.
Then, an iterate of $f$ is leaf conjugate to the time-one map
of a suspension Anosov flow.
\end{thm}

The main goal of this survey is to outline ideas in the proofs of
\cref{thm:nil,thm:sol}
and the results on which they rely.
\Cref{s.classiffoliations} details classification results for foliations
on 3-manifolds with solvable fundamental group.
\Cref{sec:leafconj} gives an
overview of the techniques to establish leaf conjugacy
which are common to all of the proofs.
The proof of \cref{thm:nil} is split into three special cases,
handled in 
\cref{sec:skew,sec:da,sec:nil}.
Finally,
\cref{thm:sol} is reviewed in \cref{sec:sol}.

\medskip

After manifolds with solvable fundamental group,
the next family of manifolds in terms of complexity are Seifert (fiber) spaces.
These generalize circle bundles over surfaces.
Ghys showed that, up to finite cover,
any Anosov flow on a circle bundle is orbit equivalent to a geodesic flow
\cite{Ghys}.
Barbot generalized this result to classify Anosov flows on all Seifert spaces
\cite{Barbot}.
In recent work with M.~Shannon, the authors established
the following \cite{HPS-New}.

\begin{thm} \label{thm:phseifert}
Let $M$ be a Seifert space whose fundamental group has exponential growth.
Then, $M$ admits a transitive partially hyperbolic diffeomorphism if and
only if it admits an Anosov flow.
\end{thm}

\begin{cor} \label{cor:noproduct}
If $\Sigma$ is a surface of genus $g \ge 2$, then $\Sigma
\times S^1$ does not admit a transitive partially hyperbolic diffeomorphism.
\end{cor}

We present an idea of the proofs in
\cref{sec:seifert}.
Note that \cref{thm:phseifert} does not provide a classification up to leaf
conjugacy.
At the moment of writing, this seems to be a difficult problem related to the recent discovery of new
examples on Seifert spaces which are not isotopic to the identity
\cite{BGP,BGHP}.
We survey these examples in \cref{ss-examples} along with other examples
\cite{BPP} which must be taken into account 
if one hopes to establish a classification
up to leaf conjugacy in full generality.  

Other classification results are presented in Section \ref{sec:others}. 

As we present the main ideas behind these results and
their proof, we also state several questions to help readers grasp the
ongoing nature of this program.
Hopefully, it serves to interest readers
and prompts them to join in this exciting quest.

\section{Classification of Anosov diffeomorphisms} \label{sec:classanosov} %

A key tool in classification results for both hyperbolic and partially
hyperbolic systems is the Franks semiconjugacy presented in this section.
We also briefly describe the classification problem for Anosov
diffeomorphisms.  

\begin{thm}
    [Franks] \label{thm:franks}
    Suppose $A:\bbT^d \to \bbT^d$ is a linear hyperbolic toral automorphism
    and $f:M \to M$ and $h_0:M \to \bbT^d$ are continuous maps
    such that
    \[
            \begin{simpleCD}
            {\pi_1(M)}
            {f_*}
            {\pi_1(M)}
            {h_{0*}}
            {}
            {h_{0*}}
            {\pi_1(\bbT^d)}
            {A_*}
            {\pi_1(\bbT^d)}  \end{simpleCD}
    \]
    commutes.  Then, there is a continuous map $h:M \to \bbT^d$
    homotopic to $h_0$ such that
    \[
            \begin{simpleCD}
            {M}
            {f}
            {M}
            {h}
            {}
            {h}
            {\bbT^d}
            {A}
            {\bbT^d}  \end{simpleCD}
    \]
    commutes.
    If $h_{0*}$ is non-zero, then $h$ is surjective.
\end{thm}

Since ideas arising in the proof are fundamental to later results for
partially hyperbolic systems,
we give a sketch of the proof.
Throughout this survey, we present ``ideas of proofs''
in order to give the reader an intuition as to why a result might hold.
These should not be regarded as rigorous proofs,
and the reader is most welcome to skip over them.

\begin{idea}[of \cref{thm:franks}]
    For simplicity,
    we only consider the case $M = \bbT^d = \bbT^2$ and where
    $h_0:\bbT^2 \to \bbT^2$
    is the identity map.
    Then,
    $f:\bbT^2 \to \bbT^2$ lifts to a map $\tf:\bbR^2 \to \bbR^2$
    which is a finite distance from the lifted linear map
    $A:\bbR^2 \to \bbR^2$.
    As $A$ is hyperbolic,
    there is a linear unstable foliation $\tAu$ of $\bbR^2$
    where each leaf $L \in \tAu$ is a translate of
    the unstable eigenspace of $A$.
    The set
    \[
        Q^u = \{ h:\tM \to \tAu: \sup_{x \in \bbR^2} \dist(x,h(x)) < \infty \}
    \]
    is a complete metric space when equipped with the metric
    \[
        D(h_1, h_2) = \sup_{x \in \bbR^2} \dist(h_1(x), h_2(x)).
    \]
    As $\tf$ and $A$ are at finite distance,
    the function $F:Q^u \to Q^u,\ \ h \mapsto A \circ h \circ \tf \inv$
    is well defined.
    If $\lam < 1$ is the stable eigenvalue of $A$,
    then
    \[    
        \dist(A(\cL_1), A(\cL_2)) = \lam \dist(\cL_1, \cL_2)
    \]
    for any two (linear) leaves $\cL_1, \cL_2 \in \tAu$
    and thus
    \[    
        D(F(h_1), F(h_2)) = \lam \, D(h_1,h_2)
    \]
    for any $h_1,h_2 \in Q^u$.
    As $F$ is a contraction, it has a unique fixed point $H^u \in Q^u$.
    This function $H^u:\tM \to \tAu$
    assigns to each point $x$ the unique leaf $\cL$
    such that
    $\dist(\tf^n(x), A^n(\cL))$
    is bounded for all $n \in \bbZ$.

    Note that if $X$ is a closed, non-empty subset of $Q^u$
    and $F(X)=X$,
    it follows that $H^u \in X$.
    From this one can show that $H^u$ is continuous
    and commutes with every deck transformation $\alpha \in \pi_1(\bbT^2)$.

    An analogous function $H^s:\tM \to \tAs$ exists.
    Define $H(x) = H^u(x) \cap H^s(x)$.
    Then $H:\bbR^2 \to \bbR^2$ satisfies $H \circ \tf = A \circ H$
    and quotients down to a semiconjugacy on $\bbT^2$.
\end{idea}
The above theorem may be used to find a semiconjugacy $h$
from a non-linear Anosov systems $f$ to its linear part $A$.
To establish that $h$ is a true topological conjugacy between $f$ and $A$,
one must also show that $h$ is invertible.
Franks showed that under certain assumptions for $f$, the above proof may be
adapted to construct a semiconjugacy in the opposite direction.
That is, to find $k:\bbT^d \to M$ where $k A = f k$ and
then $h \inv = k$.
He used this to give a classification of Anosov diffeomorphisms on surfaces \cite{FranksAnosov}.
He also gave specific conditions for classification
that were later verified by S.~Newhouse \cite{NewAnosov} and A.~Manning \cite{Manning}
to establish the following results.

\begin{thm}
    [Franks--Newhouse]
    If $f:M \to M$ is an Anosov diffeomorphism with
    $\dim \Eu = 1$, then $M$ is homeomorphic to a torus $\bbT^d$.
\end{thm}
\begin{thm}
    [Franks--Manning] \label{thm:franksmanning}
    If $f:\bbT^d \to \bbT^d$ is an Anosov diffeomorphism,
    then $f$ is topologically conjugate to a
    hyperbolic toral automorphism $A:\bbT^d \to \bbT^d$.
\end{thm}
The Franks--Manning result also holds for Anosov systems supported on
nilmanifolds and infranilmanifolds showing they are conjugate to simple
algebraic examples.
The fundamental question is if there are any other possibilities.

\begin{quest}
    Is every Anosov diffeomorphism topologically conjugate
    to an algebraic map on a infranilmanifold?
\end{quest}  %
%
There do exist examples on manifolds
which are homeomorphic but not \emph{diffeomorphic} to tori or infranilmanifolds \cite{FarrellJones,FarrellGogolev-exotic}. 

By putting conformality or ``pinching'' conditions on the restrictions
of the derivative $Df:TM \to TM$ to $\Eu$ and $\Es$,
researchers have established classification results in several settings (see for example
 \cite{BrinManning,Sadovskaya}).
 
The proof of \cref{thm:franksmanning} considers the number of periodic
points of a given period and relates this to the action of $f$ on homology
groups $H_k(M, \bbR)$ using the Lefschetz fixed point theorem.
There has been some success on applying this technique for systems defined on
other manifolds.
For an $n$-dimensional sphere,
the dimension of a homology group $H_k(\bbS^n, \bbR)$ is either zero or one
dimensional, and this may be used to show that no sphere supports an Anosov
diffeomorphism.
Using further techniques from algebraic topology,
A.~Gogolev and F.~Rodriguez Hertz have ruled out Anosov diffeomorphisms
on many families of manifolds, including certain sphere bundles
and products of spheres with tori \cite{GRH} (see also \cite{GL} and references therein for more known obstructions). %

In the transitive case, other obstructions are obtained by studying the way the strong foliations ``wrap'' around the homology of the manifold via the study of \emph{currents} (see \cite{RuelleSullivan,GRH}). 

The following important question remains open.

\begin{quest} 
Is there a simply connected compact manifold admitting an Anosov
diffeomorphism? 
\end{quest}

Spheres do not admit Anosov diffeomorphisms because their homology is very
small, but some of their products are already difficult to understand. In $\bbS^2\times
\bbS^2$ the proof that there are no Anosov diffeomorphisms is very recent
\cite{GRH}. It is unknown if $\bbS^3 \times \bbS^3$ admits Anosov diffeomorphisms.

\section{Integrability of the invariant bundles}\label{sec:integrability} %

This section discusses the work of M.~Brin, D.~Burago, and S.~Ivanov
on the existence of branching foliations
for 3-dimensional partially hyperbolic systems. %

\subsection{Unique integrability of strong bundles} 

As with Anosov diffeomorphisms, invariant foliations play a fundamental role in the study of partially hyperbolic dynamics.
We first present what could now be considered a classical result for partially
hyperbolic diffeomorphisms \cite{HPS}.

\begin{thm}[Strong Foliations]\label{Teo-StrongUnstableManifold} Let $f: M \to M$ be a partially hyperbolic diffeomorphism. Then, there exist $f$-invariant foliations $\cW^s$ and $\cW^u$ tangent to the bundles $E^s$ and $E^u$. Moreover, the bundles $E^s$ and $E^u$ are uniquely integrable. 
\end{thm}

\begin{idea}
If the stable bundle is not uniquely integrable at a point $p$,
    then there are two distinct paths $\alpha$ and $\beta$ starting at $p$
    and tangent to $\Es$.
    Pick distinct points $x$ on $\alpha$ and $y$ on $\beta$ near $p$.
    As the stable bundle is contracting, $f^k(x)$ is close to $f^k(p)$
    for all positive $k$ and there is a vector $\hat x_k \in T_{f^k(p)}M$
    such that the exponential map based
    at the point $f^k(p)$ satisfies
    $\exp(\hat x_k) = f^k(x)$.
    A sequence $\hat y_k$ may be defined analogously
    and the same properties hold.

    We may assume that $x$ and $y$ were chosen so that 
    the angle that $\hat x_0 - \hat y_0$ makes with $\Ecu(p)$ is very small.
    Since $Df$ contracts vectors in $\Es$ more sharply than those in $\Ecu$
    and since $\hat x_{k+1} - \hat y_{k+1}$ is approximately equal to
    $Df(\hat x_k - \hat y_k)$, one can show that the angle that
    $\hat x_k - \hat y_k$ makes with $\Ecu(f^k(p))$ converges to zero
    as $k \to \infty$.
    The sequences $\|\hat x_k\|$ and $\|\hat y_k\|$ each tend to zero
    at a rate associated to the stable direction $\Es(f^k(p))$.
    This is strictly faster than the rate of contraction of
    $\|\hat x_k - \hat y_k\|$ which is associated to $\Ecu(f^k(p))$.
    This gives a contradiction.
    \begin{figure}
        [t]
        \begin{center}
            \includegraphics{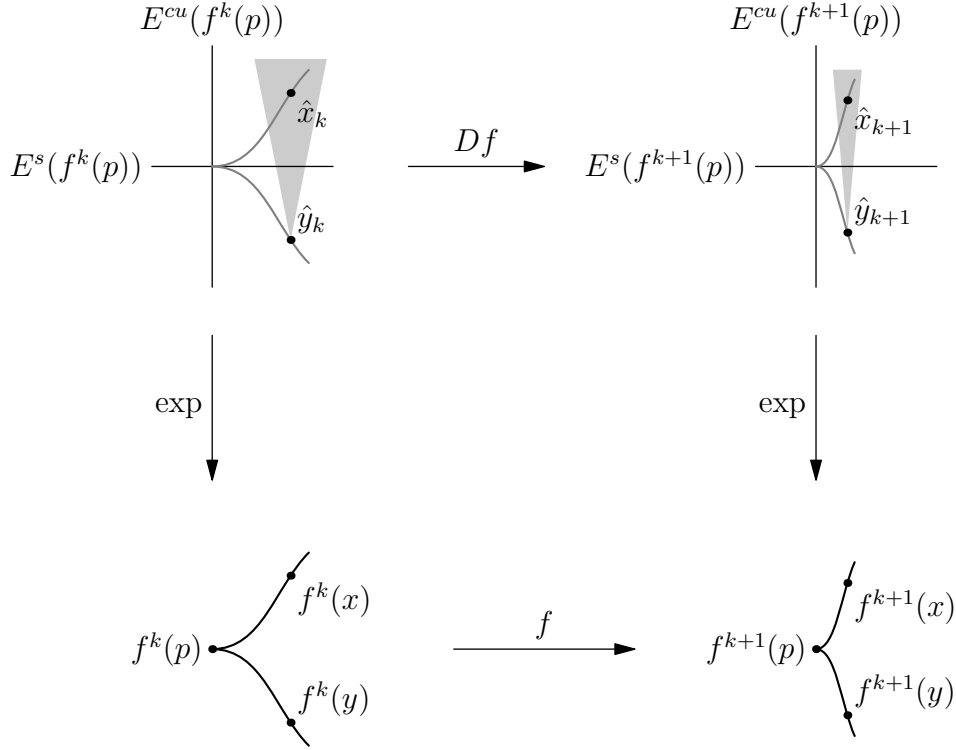}  \end{center}
        \caption{
        A depiction of the proof of \cref{Teo-StrongUnstableManifold}.
        While this diagram does not commute, it comes arbitrarily close to 
        commuting as $k \to \infty$ and the distances between the points shrink to zero.
        The two triangles in the top left and top
        right represent small cones containing the $\Ecu$ direction.
        }
        \label{fig:exp}
    \end{figure}
\end{idea}

An important simple remark is the following: 

\begin{prop} \label{prop:wulines}
    The foliations $\Ws$ and $\Wu$ have no compact leaves.
\end{prop}
\begin{idea} We argue for $\Ws$. The argument for $\Wu$ is symmetric. 
    As in all of this section, we assume $\dim \Es = 1$.
    Suppose a compact leaf $L$ exists.
    Then $L$ is a $C^1$ circle of some finite length.
    As $Df$ contracts stable vectors,
    the length of $f^{n}(L)$ shrinks to zero
    as $n \to \infty$.
    For large enough $n$,
    $f^{n}(L)$ is a circle lying inside a single foliation
    chart of $\Ws$, a contradiction.
\end{idea}

\subsection{Branching foliations} 

The main result of \cite{BI} may be summarized as follows: 

\begin{thm} \label{thm:BI}
    Let $f:M \to M$ be a partially hyperbolic diffeomorphism
    on a compact 3-manifold such that the bundles $\Eu$, $\Ec$, and $\Es$
    are oriented and $f$ preserves these orientations.
    Then,
    there is a collection $\Fcs$ of immersed surfaces tangent to $\Ecs$
    such that every $x \in M$ lies in at least one surface and
    no two surfaces in $\Fcs$ topologically cross.
    The collection is $f$-invariant{:}
    $L \in \Fcs$ if and only if $f(L) \in \Fcs$.
\end{thm}

\begin{figure}
    [t]
    \begin{center}
        \includegraphics{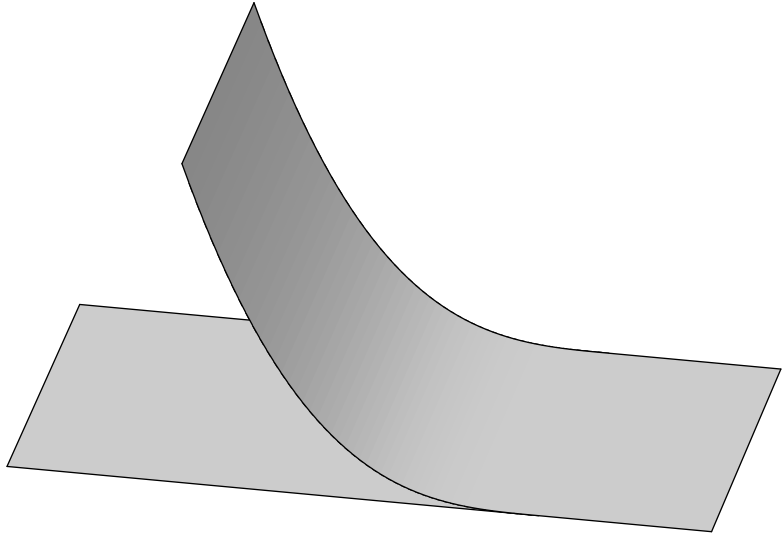}  \end{center}
    \caption{
    Merging leaves in a branching foliation.
    }
    \label{fig:bran3}
\end{figure}

The collection of surfaces $\Fcs$ in \cref{thm:BI}
is called a
\emph{branching foliation}
and its elements are called leaves.
Distinct leaves may intersect
without crossing, as depicted in
\cref{fig:bran3}.
If $f$ satisfies the hypotheses of \cref{thm:BI}
then so does $f \inv$, and this means there is also a branching foliation
$\Fcu$ tangent to $\Ecu$.
The conditions on orientation in the theorem do not hold in
general.
However, they may always be achieved by lifting an iterate of $f$ to a finite
cover.
Throughout this survey paper, we lift to a finite cover
and replace $f$ by an iterate $f^n$ whenever convenient.
Thus, for the remainder of the paper we always assume that the orientation
conditions in \cref{thm:BI} are satisfied.

To be able to apply foliation theory, the following result from \cite{BI} is also important. 

\begin{thm} \label{thm:separateleaves}
    Let $\ep > 0$ and let $\Fcs$ be a branching foliation as in
    \cref{thm:BI}.
    Then there is a true foliation $\Fcsep$ and a continuous surjective
    map $h_\ep:M \to M$ with the following properties{:}
    \pagebreak[3]
    \begin{enumerate}
        \item every leaf of $\Fcs$ is the image $h_\ep(L)$ of a leaf
        $L \in \Fcsep$,

        \item the tangent plane $T_x \Fcsep$ is $\ep$-close to $\Ecs$
        for all $x \in M$,
        and
        
        \item %
        $d(h_\ep(x), x) < \ep$ for all $x \in M$.
      \end{enumerate}  \end{thm}

For the interested reader,
we spend the remainder of this section
giving a rough outline of the ideas of the proofs of
\cref{thm:BI,thm:separateleaves}.

\subsection{Some ideas of the proof of Theorem \ref{thm:BI}} 

In the three-dimensional setting,
the center bundle is one-dimensional and therefore Peano's existence result applies. This makes one wonder if it is possible to use the unique integrability of $\Es$ to obtain a foliation tangent to $\Ecs$. However, even if $\Ec$ were uniquely integrable, it may well be that $\Ecs$ is not integrable. This is the case for example for $\Es \oplus \Eu$, which is rarely integrable even though each of $\Es$ and $\Eu$ is uniquely integrable.

Using dynamical arguments, it is possible to show that a topological type of Frobenius argument holds for the bundle $\Ecs$.
This is related to the fact that $Df\inv$ uniformly and strongly contracts
vectors transverse to $\Ecs$.

\begin{prop}[Proposition 3.1 of \cite{BI}] \label{prop:tansurf}
    Let $\gam$ be a curve tangent to $\Ecs$ and nowhere tangent to $\Es$.
    Then, the set
    $\bigcup_{x \in \gam} \Ws(x)$
    is an immersed $C^1$ surface tangent to $\Ecs$.
\end{prop}
\begin{idea}%
    This is a graph transform argument.
    For $x \in M$, let $\Ws_{\loc}(x)$ denote the stable curve of length one
    centered at $x$.
    For each integer $n > 0$,
    let $S_n$ be a $C^1$ surface through $f^n(\gam)$ which is approximately
    tangent to $\Ecs$ and which is $C^0$ close
    to $\bigcup_{x \in f^n(\gam)} \Ws_{\loc}(x)$.

    As $Df \inv$ strongly contracts vectors in $\Eu$,
    tangent planes to the surfaces $f^{-n}(S_n)$
    converge uniformly to $\Ecs$ as $n \to \infty$.
    Under the appropriate topology,
    the surfaces $f^{-n}(S_n)$ themselves
    converge to a surface tangent to $\Ecs$
    and one can show that this limit surface is equal to
    $\bigcup_{x \in \gam} \Ws(x)$.
\end{idea}

The main technical core of the paper \cite{BI} deals with the construction of the collection of complete surfaces by making use of the previous proposition. %

Using Proposition \ref{prop:tansurf} we know that given a center curve, one can saturate by strong stable manifolds in order to obtain a surface tangent to $\Ecs$. Notice that the size of such a surface for a given size of center curve is uniform. By this we mean that we can construct a uniform sized patch around each point of the center curve (given by the fact that the bundles have uniformly bounded angles). 

\pagebreak[2] %

The main difficulties are the following:

\begin{enumerate}
 \item one cannot ensure that the saturation of a (complete) center curve is complete with the induced topology, 
 \item if one considers all center curves one would get an $f$-invariant collection of surfaces, but they might topologically cross; alternatively, one could consider only one center curve at each point so that the elements of the collection do not cross, but $f$-invariance becomes more difficult to establish. 
\end{enumerate}

The proof uses the orientation of $\Eu$ and the fact that $Df$ preserves this orientation. This allows one to choose, for each point, its lowermost  local surface tangent to $\Ecs$ and one will get both non-crossing and $f$-invariance.  Let us elaborate a bit on this.
    Suppose $\sig$ is a small disk transverse to $\Es$ (and thus also transverse to $\Ecs$)
    and with $p$ in the interior of the disk.
    Let $E_\sig$ denote the intersection of $\Ecs$ with
    the tangent bundle $T \sig$.
    This intersection is a line field on the disk,
    and the orientation of $\Ec$ gives an orientation to $E_\sig$.
    With respect to this orientation,
    we may consider
    \emph{forward} and \emph{backward} curves tangent to $E_\sig$.
    The orientation of $\Eu$ gives a transverse orientation
    to $\Ecs$ which in turn gives a transverse orientation to $E_\sig$.
    Viewing $E_\sig$ as roughly horizontal,
    the transverse orientation
    gives a notion of ``up'' and ``down'' on the disk,
    and we may consider the \emph{lowest} forward curve $\gam$
    starting at $p$ and tangent to $E_\sig$.
    By \cref{prop:tansurf},
    $S := \bigcup_{x \in \gam} \Ws(x)$
    is an immersed surface tangent to $\Ecs$.
    This surface may not be complete and
    it has a boundary consisting of stable leaves.

    Consider now a point $\hat p \in \partial S$
    and a small immersed disk $\hat \sig$ containing $\hat p$
    and transverse to $\Es$.
    Similar to before, define an oriented line field $E_{\hat \sig}$
    and let $\hat \gam$ be the lowest forward curve tangent $E_{\hat \sig}$.
    If $\hat \gam \subset S$, there is nothing to do at $\hat p$.
    Otherwise, define 
    $\hat S := \bigcup_{x \in \hat \gam} \Ws(x)$
    and glue $S$ to $\hat S$ in order construct a larger surface
    tangent to $\Ecs$.

    Continuing in this way,
    one builds a surface $S_{p,0}$ tangent to $\Ecs$
    such that $\partial S_{p,0}$ consists of stable curves
    and the orientation of $\Ec$ points into $S_{p,0}$
    at all points on the boundary.
    The surface $S_{p,0}$
    depends only on the point $p$ and the oriented bundles
    $\Ec$, $\Es$, and $\Eu$.
    In particular,
    it does not depend on the various choices of disks and points $\sig$,
    $\hat p$, $\hat \sig$, et cetera that occur in the construction.
    Since $\Ec$ points into $S_{p,0}$,
    there is a uniform lower bound on the distance between points
    lying on distinct components of the boundary.

    If $p,q \in M$,
    then the surfaces $S_{p,0}$ and $S_{q,0}$ do not topologically cross.
    This follows from the method of construction,
    though the proof is not immediate.
    Since the partially hyperbolic diffeomorphism $f$ preserves
    the subbundles 
    $\Ec$, $\Es$, and $\Eu$
    and their orientations,
    it follows that $f(S_{p,0}) = S_{f(p),0}.$
    \begin{figure}
        [p!]
        \begin{center}
            \includegraphics{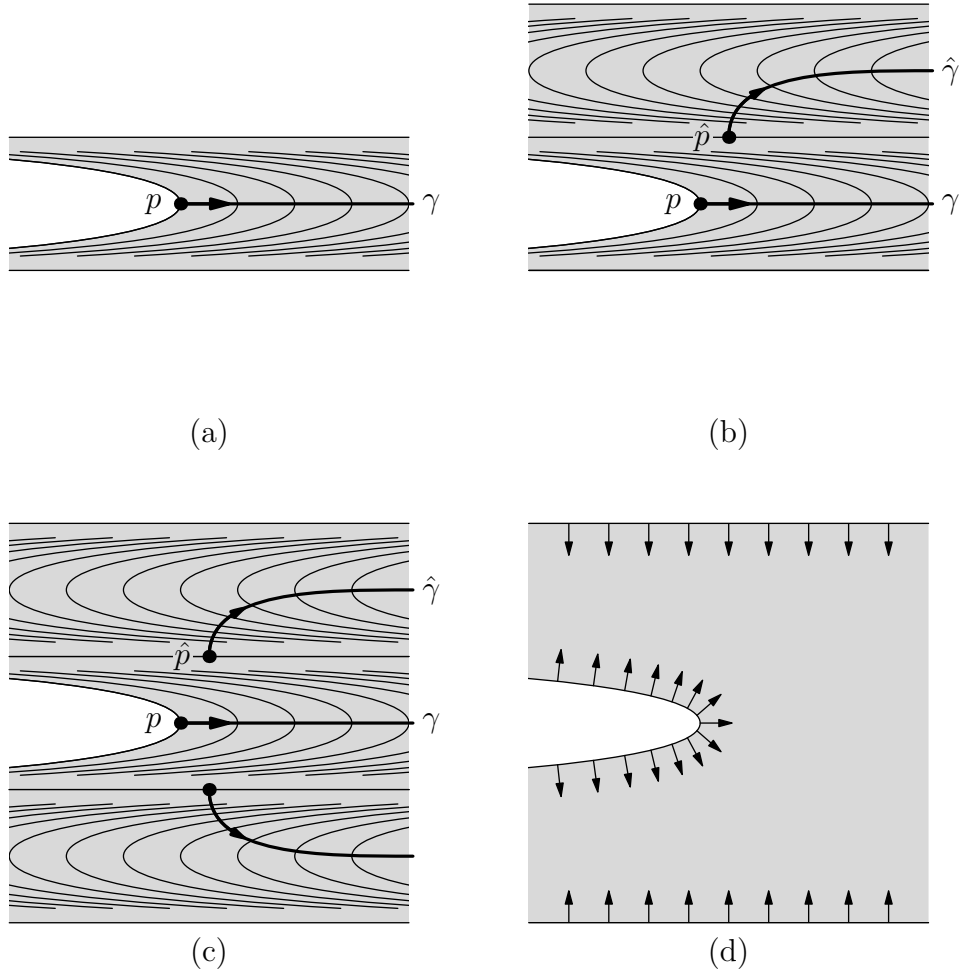}  \end{center}
        \caption{
        A depiction of a center-stable leaf being constructed in stages.
        Starting at a point $p$, take a center curve $\gamma$ which is locally
        the lowest possible forward center curve through each of its points.
        Let $S$ be the union of
        stable curves passing through $\gamma$ (a).
        At a point $\hat p$ on the boundary of $S$, the center bundle points
        away from $S$ and so one considers a lowest forward curve 
        $\hat \gamma$ extending from $\hat p$.  The stable curves through
        $\hat \gamma$ give a surface $\hat S$ which is glued to $S$
        (b).
        By the same process, a similar surface is glued to the other
        boundary component of $S$ (c).
        Now in this example, every point on the boundary of the constructed
        surface $S_{p,0}$ points into the surface and this stage of the
        construction is finished (d).
        The next stage will attach surfaces by considering
        center curves with the opposite orientation.
        }
        \label{fig:buildsurf}
    \end{figure}

    Now consider $p$ and $\sig$ as before and let $\gam_1$ be a backward
    curve starting at $p$ and tangent to $E_\sig$.
    Further, assume $\gam_1$ is the highest possible curve
    which does not topologically cross $S_{q,0}$ for any $q \in M$.
    Attach the surface
    $\bigcup_{x \in \gam_1} \Ws(x)$
    to $S_{p,0}$.
    Continuing in this way by attaching the stable saturates of
    highest backward curves,
    one builds a surface $S_{p,1} \supset S_{p,0}$
    so that the orientation of $\Ec$
    now points away from $S_{p,1}$ at all points on the boundary.
    If $p,q \in M$, then $S_{p,1}$ and $S_{q,1}$
    do not topologically cross.

    By alternating between forward and backward curves
    and at each step taking care not to cross a surface constructed in a
    previous step, one builds a nested sequence of surfaces $S_{p,k}$
    which in the limit as $k \to \infty$
    yields a complete immersed surface through the point $p$.
    We refer the reader to \cite[Sections 4--6]{BI} for (many) more details. 

\medskip

Finally, let us mention that it is proven in \cite[Lemma 7.1]{BI} that the collection can be chosen to satisfy that if $x_n \to x$ and $L_n$ are surfaces of the collection through $x_n$, then the sequence $L_n$ converges to a surface $L$ in the collection. 

\subsection{Pulling leaves apart}

We give here an indication of the proof of Theorem \ref{thm:separateleaves}. 
 The idea is to flow leaves by a small amount depending on the point and the leaf in order to ``separate them''. Of course, it is essential that leaves do not cross in order that this will work. 

This can be done without problem in the core of a coordinate chart where
    we may consider a branching foliation on the unit cube
    $[0,1]^3 \subset \bbR^3$
    where each leaf is the graph of a function $[0,1]^2 \to \bbR$.
    The condition that leaves do not topologically cross
    implies that if two leaves are given by the graphs
    of $\phi$ and $\psi$, then
    one of the inequalities
    $\phi  \le  \psi$ or $\psi  \le  \phi$
    holds pointwise for all $(x,y) \in [0,1]^2$.
    Thus, the leaves are a totally ordered set.
    Properties of ordered sets imply that we may associate to each leaf
    a number $t \in [0,1]$ such that if $\graph(\phi_t)$ is the leaf associated
    to $t$, then
    \[
        \phi_s  \le  \phi_t
        \quad \text{if and only if} \quad
        s  \le  t.
    \]
    Now take $\delta > 0$ and define functions
    \begin{math}
        \hat \phi_t(x,y) = \phi_t(x,y) + \delta t.  \end{math}
    This family of functions has the stronger property that
    \[
        \hat \phi_s < \hat \phi_t
        \quad \text{if and only if} \quad
        s < t.
    \]
    The resulting graphs are therefore pairwise disjoint and form a true
    foliation.     This solves the problem locally.

    To solve the problem globally,
    cover $M$ with a finite number of foliation charts and perform
    this procedure in each chart in parallel, using a partition of unity
    to make the result continuous and taking care to ensure that changes in
    one chart are small enough that they do not cause problems in other
    charts. See \cite[Section 7]{BI} for details.

The above sketch implicitly used the following property of branching
foliations.

\begin{prop} \label{prop:nobranch}
    If the leaves of a branching foliation
    are pairwise disjoint, then the branching foliation
    is a true foliation.
\end{prop}
See the discussion in \cite[Remark 1.10]{BW} for further details.
For the classification results in this survey,
this property will allow us to establish dynamical coherence 
as a step towards leaf conjugacy.

\subsection{Dynamical coherence} 

Integrability of $\Ecs=\Es \oplus \Ec $ and $\Ecu=\Ec \oplus \Eu$ into $f$-invariant foliations is called \emph{dynamical coherence} as introduced in the first section of this survey (see also \cite{BuW2}).
It is unknown in general if transitive partially hyperbolic diffeomorphisms in dimension 3 are dynamically coherent.
There exist non-transitive examples of non-dynamically coherent partially hyperbolic diffeomorphisms \cite{HHU}.

We refer the reader to \cite{CHHU} for a survey on this problem in
3-dimensional manifolds. In this survey, integrability of the bundles appears as a step in the quest for showing leaf conjugacy to the known models.

As a consequence of Proposition \ref{prop:nobranch} we have:

\begin{cor}
    For a partially hyperbolic diffeomorphism $f$ on a compact 3-manifold,
    if the leaves of the branching foliation $\Fcs$ are pairwise disjoint and
    the leaves of  the branching foliation $\Fcu$ are pairwise disjoint,
    then $f$ is dynamically coherent.
\end{cor}
\bigskip{}

\section{Reebless foliations} \label{sec:reebless} %

In the case of a strongly partially hyperbolic system,
the tangent space $TM$ splits into three subbundles,
each of dimension at least one.
The manifold $M$ must therefore have dimension at least three.
For a long time,
this was the only known obstruction;
it was an open question whether every manifold $M$
with $\dim M  \ge  3$
supported a partially hyperbolic diffeomorphism.

This was finally answered by Brin, Burago, and Ivanov,
who used the branching foliations described in \cref{sec:integrability}
to show the following.

\begin{thm}
    There are no partially hyperbolic diffeomorphisms on the 3-sphere.
\end{thm}
The proof relies on the theory of compact leaves and Reeb components
developed by Novikov \cite{novikov1965}.
We give a brief overview of this theory
and refer the reader to \cite{CandelConlon} for a detailed treatment.

\begin{figure}
    [t]
    \begin{center}
        \includegraphics{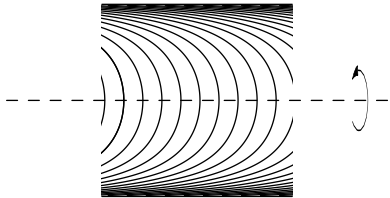}  \end{center}
    \caption{
    A construction of a Reeb component.
    }
    \label{fig:reeba}
      \end{figure}

First,
consider an annulus $\bbS^1 \times [0,1]$
constructed by identifying the left and right sides of a square.
On this annulus, consider a foliation as in \cref{fig:reeba}
with two circle leaves on the boundary and
where all other leaves are lines.
Finally,
rotate $\bbS^1 \times [0,1]$ around the axis
$\bbS^1 \times \tfrac{1}{2}$ to
construct a solid torus $\bbS^1 \times D^2$.
The resulting two dimensional foliation
has the boundary as a compact leaf
and all other leaves are topological planes.
The solid torus equipped with this foliation is called a
Reeb component.

Suppose $M$ is a 3-manifold and $\cF$ a two-dimensional foliation on $M$.
If there is a solid torus $R \subset M$
such that the foliation $\cF$ restricted
to $R$ is equivalent to the construction above,
then $R$ is a Reeb component for $\cF$.
The embedding $i: \partial R \hookrightarrow M$
induces a homomorphism $i_*: \pi_1(\partial R) \to \pi_1(M)$.
Since the embedding factors as a composition
$\partial R \hookrightarrow R \hookrightarrow M$,
it is easy to see that $i_*$
is not injective.
Novikov proved the much more difficult converse.

\begin{thm}
    \label{thm:piinj}
    Suppose $\cF$ is a two-dimensional foliation of a 3-manifold.
    If there is a leaf $L$ such that $\pi_1(L) \to \pi_1(M)$
    is not injective,
    then $\cF$ has a Reeb component.
\end{thm}
This was originally proved for smooth foliations and later extended to the
$C^0$ case by Solodov \cite{Solodov}. 

Call a foliation
\emph{Reebless}
if it has no Reeb components.
The above results then say that a foliation is Reebless
if and only if
every leaf is $\pi_1$-injectively immersed in $M$.

\begin{prop} \label{prop:reebcircle}
    Suppose $\cF$ is a two-dimensional foliation
    and $\cW$ a transverse one-dimensional foliation
    on a 3-manifold.
    If $\cF$ has a Reeb component,
    then $\cW$ has a compact leaf.
\end{prop}

\begin{idea}
    Assume $\cW$ is orientable, so that there is a flow
    whose orbits are exactly the leaves of $\cW$.
    Suppose $\cF$ has a Reeb component $R$
    and assume on the boundary of $R$ that
    the direction of the flow is pointing into the solid torus.
    Then any orbit which enters $R$ cannot leave $R$.
    With some work, one can show that there a compact disk lying inside
    a leaf in $R$ to which every orbit of the flow returns.
    The Brouwer fixed point theorem then implies the flow has a periodic orbit.
\end{idea}

Novikov also established a form of converse to \cref{prop:reebcircle}.
Given a foliation $\cF$,
call a circle $S \subset M$ a
\emph{transverse contractible cycle}
if is it contractible (that is, homotopic in $M$ to a constant path)
and transverse to $\cF$.
(While the original $S$ is transverse to $\cF$,
the circles in the homotopy need not be.)

\begin{prop}\label{prop:transversecycle}
    A Reebless foliation has no transverse contractible cycle.
\end{prop}
For the remainder of the section,
we consider foliations for the specific case
of a partially hyperbolic diffeomorphism $f$ on a 3-manifold $M$.
Let $\Wu$ be the unstable foliation. Recall that Proposition \ref{prop:wulines} implies that $\Wu$ has no circle leaves.

By \cref{thm:BI,thm:separateleaves},
assume $\Fcs$ is a branching foliation approximated by a true foliation
$\Fcsep$.  Take $\ep$ small enough that $\Fcsep$ is transverse to $\Wu$.
\Cref{prop:reebcircle} and \cref{prop:wulines} immediately give the following.

\begin{prop}
    The foliation $\Fcsep$ is Reebless.
\end{prop}
Let $\tM$ be the universal cover of $M$.
As $\Fcsep$ is a true foliation,
lifting it to a foliation $\tFcsep$ on $\tM$ is straightforward.
To lift the branching foliation,
let $h_\ep:M \to M$ be the map which takes leaves of $\Fcsep$
to leaves of $\Fcs$.
Lift $h_\ep$ to $\tilh_\ep: \tM \to \tM$.
Then the leaves of $\tFcs$ are exactly the surfaces of the form
$\tilh_\ep(\tL)$ where $\tL \in \tFcsep$.

\begin{prop}
    Every leaf of $\tFcsep$ is a plane.
\end{prop}
\begin{idea}
    Suppose $\tL$ is a leaf of $\tFcsep$
    and $\tgam$ is a closed path in $\tL$.
    Let $L$ and $\gam$ be their projections in $M$.
    Since $\tgam$ is contractible in $\tM$,
    $\gam$ is contractible in $M$.
    By \cref{thm:piinj},
    $\gam$ is also contractible in $L$.
    The homotopy in $L$ taking $\gam$ to a point
    lifts to a homotopy in $\tL$ taking $\tgam$ to a point.
    This shows that $\tL$ is simply connected.

    The surface $\tL$ is then either a plane or a 2-sphere.
    The map $\tilh_\ep$ restricts to a homeomorphism between $\tL$
    and a leaf of $\Fcsep$.
    This leaf $\tilh_\ep(\tL)$ has a splitting $\Ec \oplus \Es$ of its tangent
    bundle and therefore cannot be a sphere.
\end{idea}
Palmeira \cite{palmeira} showed that if a simply connected $d$-manifold $\tM$ is foliated
by planes of dimension $d-1$, then $\tM$ is homeomorphic to $\bbR^d$. 
Thus in the current setting, $\tM$ must be homeomorphic to $\bbR^3$.
Due to its importance, we state this as a self-contained theorem.

\begin{thm}
    If a 3-manifold supports a partially hyperbolic diffeomorphism,
    its universal cover is homeomorphic to $\bbR^3$.
\end{thm}
Though $\Fcs$ may have branching,
in many ways it behaves like a Reebless foliation.

\begin{prop}
    There is no contractible cycle transverse to $\Fcs$.
\end{prop}
\begin{idea}
    Suppose $\gam$ is a transverse cycle.
    The angle between $\gam$ and $\Fcs$ is bounded away from zero,
    say by $\delta > 0$.
    Take the approximation $\Fcsep$ with $\ep < \delta$.
    Then $\gam$ is also transverse to $\Fcs$ and therefore not
    contractible.
\end{idea}
This implies a uniqueness property for intersections on the universal cover.

\begin{prop} \label{prop:uscunique}
    A leaf of $\tFcs$ intersects a leaf of $\tWu$
    in at most one point.
\end{prop}
\begin{idea}
    Suppose $x \in L \in \tFcs$
    and that $\tWu(x)$ intersects $L$ at a point $y  \ne  x$.
    A new path can be constructed by following $\tWu(x)$ from $x$
    to a point $z$ just before $y$,
    and then connecting $z$ to $x$ by a path which stays close
    to $L$ and which is everywhere transverse to $L$.
    \begin{figure}
        [t]
        \begin{center}
            \includegraphics{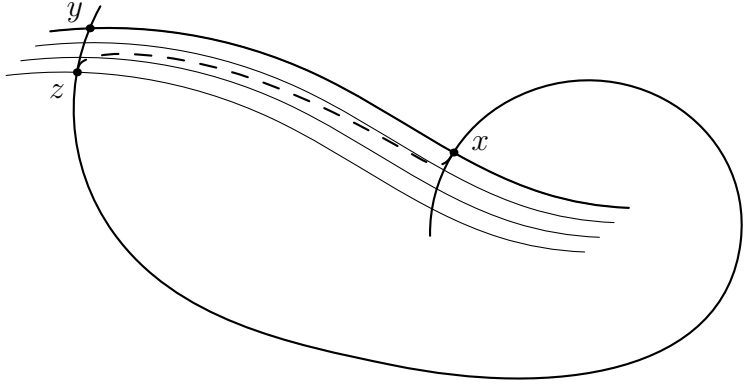}  \end{center}
        \caption{
        If an unstable leaf intersects a leaf of $\tFcs$ in distinct
        points $x$ and $y$, a transverse cycle may be constructed.
        }
        \label{fig:transverse}
          \end{figure}
    See \cref{fig:transverse}.
    The result is a cycle in $\tM$ transverse to $\tFcs$.
    This projects down to a transverse contractible cycle in $M$.
\end{idea}

\begin{cor} \label{cor:halfspace}
    Each leaf $L$ of $\tFcs$ is a properly embedded plane which
    splits $\tM$ into two half spaces.
\end{cor}
\begin{idea}
    If $L$ is not properly embedded, it accumulates at a point $x \in \tM$
    and intersects the unstable leaf through $x$ in infinitely many points.
\end{idea}

\Cref{prop:uscunique} implies that different segments of the same unstable
leaf must stay a certain distance away from each other.
This notion may be stated quantitatively.
For a subset $X \subset \tM$, let $U_1(X)$ be the neighborhood of all points
at distance at most $1$ from $X$.  That is, $y \in U_1(X)$
if $\inf_{x \in X} d(x,y) < 1$.

\begin{lem} \label{lem:volume}
    There is $C > 0$ such that if $J \subset \tM$ is an unstable curve,
    then
    \[    
        \volume(U_1(J)) > C \length(J).
    \]  \end{lem}
\begin{idea}
    For an unstable segment $J$,
    define $U^{cs}_1(J)$ as those points $x \in U_1(J)$
    such that $L \cap J  \ne  \varnothing$ for every leaf $L \in \tFcs$
    with $x \in L$.
    Since $\Ecs$ and $\Eu$ are transverse,
    one can find $\delta > 0$ such that
    $\volume U^{cs}_1(J) > \delta$ for all unstable curves $J$
    of length exactly one.

    If $\length(J) > n$
    for some integer $n$, there are pairwise disjoint subcurves
    $J_1, \ldots, J_n$
    each of length one.
    By \cref{prop:uscunique},
    the sets
    $U^{cs}_1(J_1), \ldots, U^{cs}_1(J_n)$
    are pairwise disjoint
    which implies that $\volume U_1(J) > \delta n$.
\end{idea}
If $J \subset \tM$ is an unstable curve,
then the length $\tf^n(J)$ grows exponentially fast.
If the universal cover $\tM$ is ``small''
in some geometric sense, then this is only possible if
the action $f_*$ on $\pi_1(M)$
exhibits some form of partial hyperbolicity.
To give an idea of this,
we prove a simple result on $\bbT^3$.

\begin{prop}
    If $f:\bbT^3 \to \bbT^3$ is partially hyperbolic,
    then $f_*:\pi_1(\bbT^3) \to \pi_1(\bbT^3)$
    is not the identity.
\end{prop}
\begin{idea}
    Suppose $f_*$ is the identity.
    Then, $f$ lifts to a map $\tf:\bbR^3 \to \bbR^3$
    which is a finite distance from the identity map on $\bbR^3$.
    There is $c > 0$ such that if
    $B(r) \subset \bbR^3$ denotes the ball of radius $r$,
    then
    \begin{math}
        \tf^n(B(r)) \subset B(r + c n)
    \end{math}
    for all $n  \ge  0$.
    Take an unstable segment $J \subset B(r)$.
    Then,
    \[
        \volume U_1( \tf^n(J))  \le  \volume B(r + c n + 1)
    \]
    which grows at most polynomially in $n$.
    Since $\length \tf^n(J)$ grows exponentially,
    \cref{lem:volume} gives a contradiction.
\end{idea}
Since $\pi_1(\bbT^3) \cong \bbZ^3$,
the homomorphism
$f_*$ may be regarded as a $3 \times 3$ matrix
with eigenvalues $\lam_i$ ordered by
$|\lam_1|  \le  |\lam_2|  \le  |\lam_3|$.
Then $\lim \| f^n_* \|^{\frac{1}{n}} = |\lam_3|$
by Gelfand's formula and
the above proof can be adapted to show that $|\lam_3| > 1$.
A similar study of $f \inv$ shows that $|\lam_1| < 1$.
This is explored in more detail in \cref{sec:skew,sec:da}.

\section{Interlude: weak partial hyperbolicity}\label{sec:interlude}

The arguments in the previous section all involved the one-dimensional
unstable foliation $\Wu$ and the two-dimensional branching foliation $\Fcs$.
Once we established the existence of the branching foliation
$\Fcs$, the finer splitting $\Ecs = \Ec \oplus \Es$ was not directly used.
One could therefore ask if the arguments also apply to weakly partially hyperbolic
systems.  Indeed, the following holds.

\begin{prop}\label{prop:weakfoln}
    Let $f:M \to M$ be a weakly partially hyperbolic diffeomorphism
    with a splitting of the form $TM = \Ec \oplus \Eu$
    where $\dim \Ec = 2$ and $\dim \Eu = 1$.
    Further, suppose there is a
    foliation $\cF$ transverse to $\Eu$.
    Then{:}
    \begin{enumerate}
        \item the universal cover $\tM$ is $\bbR^3$,

        \item the leaves of the lifted foliation $\tF$ are all planes, and

        \item each leaf of $\tWu$ intersects each leaf of $\tF$ at most once.
    \end{enumerate}  \end{prop}

This was used in \cite{Pot} to study dynamical coherence of weakly partially hyperbolic diffeomorphisms isotopic to a linear Anosov automorphism of $\mathbb{T}^3$ (and extended in \cite{FPS} to higher dimensions). The notion of \emph{almost dynamical coherence} was used there to refer to the existence of a foliation $\cF$ transverse to $\Eu$ as in the proposition. 

The problem is that in the weakly partially hyperbolic setting, we have no
tools at present to find such a foliation $\cF$.

\begin{quest}
    Let $f:M \to M$ be a weakly partially hyperbolic diffeomorphism
    with a splitting of the form $TM = \Ec \oplus \Eu$
    where $\dim \Ec = 2$ and $\dim \Eu = 1$.
    Is there a branching foliation tangent to $\Ec$?
    Is there a true foliation transverse to $\Eu$?
\end{quest}
It is not clear even locally if there exist small surfaces tangent to $\Ecs$. 
Some progress in this direction has been obtained by Luzzatto, Turelli and War
\cite{LuzzattoEtAl} continuing the work in \cite{Ham-Integrability} and recently improvements have been announced by Turelli and War. 

If there always exists a transverse foliation $\cF$ as in \cref{prop:weakfoln},
it would imply that any 3-manifold which supports
a robustly transitive diffeomorphism must also support a Reebless foliation.
That would definitively answer the following long standing open question.

\begin{quest}
    Does $\bbS^3$ support a robustly transitive diffeomorphism?
\end{quest}
In the measure-preserving setting, one can also ask if $\bbS^3$ supports a
stably ergodic diffeomorphism, a question which is also wide open. 

\medskip

One can consider even weaker forms of $Df$-invariant structures
and make an attempt at classification.
One such invariant structure occurs for \emph{non-uniformly hyperbolic} systems.
In this case, the splitting of the tangent bundle is only measurable instead of
continuous.
Therefore it does not impose much restrictions on the topology of the underlying manifolds, and indeed, it is shown in \cite{dolgopyatpesin} that every compact manifold admits a conservative non-uniformly hyperbolic diffeomorphism.

Another case is a \emph{dominated splitting}
where the tangent bundle has a $Df$-invariant splitting $TM =E\oplus F$
such that
\[
    ||Df^n u|| < \tfrac{1}{2} || Df^n v||
\]
for some $n \ge 1$ and all $x \in M$ and unit vectors $u \in E_x$ and $v \in F_x$.
In comparison to weak partial hyperbolicity, a dominating splitting puts no
conditions of uniform expansion or contraction on either
$E$ or $F$.
Gourmelon and Potrie are currently preparing a complete classification of the
dynamics and structure of diffeomorphisms of surfaces admitting global
dominated splittings \cite{gourmelonpotrie} which builds on the work of Pujals and Sambarino who have studied in great detail the dynamics of surface diffeomorphisms whose limit set admits a dominated splitting \cite{pujalssambarino}.

\section{Classification of foliations in certain 3-dimensional manifolds}\label{s.classiffoliations} %

\subsection{Compact leaves} 
In the partially hyperbolic context, we mostly consider invariant foliations
without compact leaves.  This is due to the following work of
Rodriguez Hertz, Rodriguez Hertz, and
Ures \cite{HHU3}.

\pagebreak[2]

\begin{thm} \label{thm:anosovtorus}
    Suppose $f$ is a partially hyperbolic diffeomorphism of a compact
    3-manifold $M$ and there is a compact surface tangent to $\Ecs$.
    Then{:}
    \begin{enumerate}
        \item there is a periodic 2-torus $T = f^k(T)$ tangent to $\Ecs$,

        \item the map $f^k|_T$ is isotopic to an Anosov diffeomorphism on $T$,
        and

        \item (up to finite cover) $M$ either is the 3-torus or
        the suspension of an Anosov diffeomorphism.
    \end{enumerate}  \end{thm}

The torus $T$ is a specific case of an \emph{Anosov torus},
a concept covered in detail in \cite{CHHU}
and so we do not focus on its properties here.

\subsection{Almost parallel foliations} 

Foliations $\cF_1$ and $\cF_2$ on a manifold $M$ are
\emph{equivalent}
if there is a homeomorphism $h:M \to M$ such that $L$ is a leaf of $\cF_1$
if and only if $h(L)$ is a leaf of $\cF_2$.
When the foliations are $C^r$ smooth for $r  \ge  2$,
equivalence is a natural notion to use.
For instance, the following holds for foliations on surfaces.

\begin{prop} \label{prop:c2foln}
    Let $\cW$ be a $C^2$ foliation without compact leaves
    defined on $\bbT^2$.
    Then, $\cW$ is equivalent to a foliation by lines of constant irrational
    slope.
\end{prop}

If $\cW$ is only $C^1$, the result does not hold, due to foliations with
exceptional minimal sets as discovered by Denjoy. %
For partially hyperbolic systems, the foliations under study
are typically $C^0$ foliations with $C^1$ leaves.

To analyse these foliations, it has turned out to be useful to consider a weaker equivalence of foliations without compact leaves \cite{Pot,HP,HP2}.
Two (branching) foliations $\cF_1$ and $\cF_2$ of $\tM$ are \emph{almost parallel} if there exists $K>0$ such that:  

\begin{itemize}
 \item for every $L_1 \in \cF_1$ there exists $L_2 \in \cF_2$ such that $d_H(L_1,L_2) \leq K$,
 \item for every $L_2 \in \cF_2$ there exists $L_1 \in \cF_1$ such that $d_H(L_1,L_2) \leq K$.
\end{itemize}

Here, $d_H$ denotes the Hausdorff distance between closed subsets of $\tM$ with a metric obtained by lifting the metric of $M$ to the universal cover. 

Note that ``almost parallel'' is an equivalence relation,
but the word ``equivalent'' is reserved for the stronger notion defined above.

\begin{prop} \label{prop:c0foln}
    Let $\cW$ be a $C^0$ foliation without compact leaves
    defined on $\bbT^2$.
    Then, $\cW$ is almost parallel to a foliation by lines of constant irrational
    slope.
\end{prop}

We leave this as an exercise to the reader.

\subsection{%
Foliations in 3-manifolds with solvable fundamental group} 

Plante classified the possible $C^2$ Reebless foliations
in 3-manifolds with solvable fundamental group \cite{Plante}.
Adapting parts of the proofs, we established
a weaker result for $C^0$ branching foliations \cite[Appendix B]{HP}.
See also \cite{BBI2,Pot,HP2}.
To state the result, we introduce the notion of a model foliation.

First, consider a linear two-dimensional foliation on $\bbR^3$.
That is, a foliation by geometric planes.
Restrict this foliation to $\bbR^2 \times [0,1]$
and quotient by the lattice
$\bbZ^2 \times \{0\}$
to produce a linear foliation $\cF$ on
$\bbT^2 \times [0,1]$.
Let $\cF_0$ and $\cF_1$ denote the restrictions of $\cF$
to the boundary components
$\bbT^2 \times \{0\}$
and
$\bbT^2 \times \{1\}$
Now consider a linear map $\varphi:\bbT^2 \to \bbT^2$
and identify $(v,1)$ with $(\varphi(v), 0)$ to produce the suspension manifold
$M_\varphi$.
If $\varphi$ maps $\cF_0$ to $\cF_1$, then the linear foliation $\cF$
yields a foliation $\cF_\varphi$ on the $M_\varphi$.
We call such a foliation $\cF_\varphi$ a
\emph{model foliation}.
In the special case that $\varphi$ is the identity,
$M=\bbT^3$ and the model foliations are precisely the linear foliations.

\begin{thm} \label{thm:models}
    Suppose that $\cF$ is a (branching) foliation on a compact 3-manifold with
    solvable fundamental group and that $\cF$ has no compact leaves.
    Then,
    (up to finite cover)
    $M$ is a torus bundle over a circle and $\cF$ is almost parallel to a model
    foliation.
\end{thm}

Some classification results also hold for Reebless foliations in other
families of 3-manifolds.
See the books \cite{CandelConlon,Calegari} for comprehensive treatments of
foliation theory in dimension 3 and beyond.
\Cref{sec:seifert} states some results for Reebless foliations on Seifert spaces
and applies these to the study of partial hyperbolicity on such manifolds.

\section{Leaf conjugacy in certain 3-manifolds}\label{sec:leafconj} %

We give a broad outline of how to establish a leaf
conjugacy such that all of the results for the 3-torus and other
nilmanifolds and solvmanifolds follow roughly the same blueprint.
The fact that these different classification results can be fit into
the same framework gives hope that some of the techniques will also apply to
other partially hyperbolic systems.

Given a partially hyperbolic diffeomorphism $f$
on a 3-manifold $M$, the first step is to
use properties of the geometry of $M$
to characterize the action $f_*:\pi_1(M) \to \pi_1(M)$ on the fundamental
group.
From this,
construct a diffeomorphism $g:M \to M$ which is
algebraic
and such that $f_*=g_*$.
Then,
use the information about $f_*$ to show that $g$ is partially hyperbolic.

Since the universal cover $\tM$ is homeomorphic to $\bbR^3$,
results in algebraic topology imply that $f$ and $g$ are homotopic
\cite[Proposition 1B.9]{hatcher}. 
Lift $f$ and $g$ to diffeomorphisms $\tf$ and $\tg$ 
which are at finite distance from each other;  
that is,
\[
    \sup_{x \in \tM} d(\tf(x), \tg(x)) < \infty.
\]

\Cref{thm:BI}
implies that there is a Reebless
branching foliation
$\Fcs$ tangent to $\Ecs$ on $M$ and invariant by $f$.
Lift $\Fcs$ to a branching foliation $\tFcs$ on $\tM$.
By \cref{thm:models},
there is a model foliation $\cA$ with lift $\tA$
such that $\tFcs$ is almost parallel to $\tA$.

Now consider a leaf $\cL \in \tA$.
It lies a finite distance
from a leaf $L \in \tFcs$.
Then $\tf(L) \in \tFcs$ lies a finite distance from $\tg(\cL)$.
Using the definition of almost parallel,
this implies that $\tg(\cL)$ is a finite distance from a leaf
of $\tA$.  In other words, $\tg(\tA)$ is almost parallel to $\tA$.
This further restricts the possibilities for $\tA$.

Using techniques specific to the manifold,
show that $\tA$ is equal to the center-stable foliation $\tAcs$
associated to $\tg$.
This establishes that $\tFcs$ is almost parallel to $\tAcs$.
Similarly, show that $\tFcu$ is almost parallel to the center-unstable
foliation $\tAcu$ of $\tg$.

Our final goal is to construct a homeomorphism $h:M \to M$ homotopic to the
identity which is a leaf conjugacy between $f$ and $g$.
Note that if such a map $h$ exists,
it lifts to a map $\tilh:\tM \to \tM$ which is a finite distance from the
identity.
Consider $p \in \tM$ and let
$\cL \in \tAcs$ be the leaf which passes through $\tilh(p)$.
The fact that $\tilh$
is a leaf conjugacy implies that $\tilh \tf^k(p) \in \tg^k(\cL)$
for all $k \in \bbZ$.
Then,
\[
    \sup_{k \in \bbZ} \dist(\tf^k(p), \tg^k(\cL)) \le
    \sup_{k \in \bbZ} d(\tf^k(p), \tilh(\tf^k(p)) <
    \infty.
\]
Since we do not have the leaf conjugacy $h$ in advance, we run this
reasoning backwards.
That is, we establish that for every $p \in \tM$ there is exactly one leaf
$\cL \in \tAcs$ which satisfies $\sup_{k \in \bbZ} \dist(\tf^k(p), \tg^k(\cL)) < \infty$.
This then defines a map $\Hcs:\tM \to \tAcs$ by setting $\Hcs(p)$ equal to $\cL$.

One must further show that the leaves of $\tFcs$ are the fibers of $\Hcs$.
Since the fibers of $\Hcs$ are disjoint, this implies that $\tFcs$ is a true
foliation.  Further, it is the unique invariant foliation, since the
definition of $\Hcs$ does not depend on the choice of $\tFcs$.

Similar reasoning gives a map $\Hcu: \tM \to \tAcu$. Assuming each
leaf of $\tAcs$ intersects each leaf of $\tAcu$ in exactly one center leaf,
these two maps combine to give $H:\tM \to \tAc$ where
$H(x) = \Hcs(x) \cap \Hcu(x)$.

Thus, each point of $x$ is mapped to a center leaf $H(x)$ of $\tg$
in such a way that $H$ is constant on center leaves of $\tf$.
The final step is to define a homeomorphism $\tilh:\tM \to \tM$ such
that $\tilh(x) \in H(x)$ for all $x \in \tM$ and in such a way
that $\tilh$ quotients down to a homeomorphism $h:M \to M$
on the original compact manifold.  This $h$ will be the desired
leaf conjugacy.

The original papers have three very different methods for constructing
$\tilh$ from $H$.  With the goal of unifying these three proofs as much as
possible, we show that the averaging technique used in \cite{HP2} can be used
in all of the cases.
This averaging method is due to Fuller \cite{fuller1965},
who used it to study cross sections of flows.
Others have employed the same technique \cite{verjovsky,Ghys,gromov}
and it appears to have been independently discovered multiple times.

The next step uses a smooth projection $P:\tM \to \bbR$ with the following
properties.
\begin{enumerate}
    \item For any $\cL \in \tAc$, $P|_\cL$ is a $C^1$ diffeomorphism.

    \item There is $T > 0$, such that if $x$ and $y$ are on the same center
    leaf $L$ of $\tf$ and their distance $d_c(x,y)$ measured along $L$ is
    exactly $T$,
    then $|P(x) - P(y)| > 1$.

    \item There is a homomorphism $\tau:\pi_1(M) \to \bbR$ such that
    \[        P \gam(x) = P(x) + \tau(x)  \]
    for all $x \in \tM$ and $\gam \in \pi_1(M)$.
\end{enumerate}

The choice of $P$ depends on the $3$-manifold $M$. 

Now consider $L \in \tWc$ and let $\alpha:\bbR \to L$ be a function which
parameterizes $L$ by arc length.
For $x \in L$, let $s = \alpha \inv(x)$ and define $\tilh(x) \in H(L)$
by the condition
\[
    P \tilh(x) =
    \frac{1}{T} \int_{0}^{T} P \alpha(s+t) \, dt.
\]
This defines a homeomorphism from $L$ to $H(L)$
which is independent of the choice of $\alpha$.
When the same construction is performed for every leaf in $\tWc$
the result is a leaf conjugacy $\tilh:\tM \to \tM$.
The existence of the homomorphism $\tau$ implies that $\tilh$ commutes with every
deck transformation, and therefore descends to a leaf conjugacy on the compact
manifold $M$.

The next four sections discuss steps to establish leaf conjugacy
which are unique to the manifolds under study.

\section{Euclidean geometry, not isotopic to Anosov} \label{sec:skew} %

Let $f:\bbT^3 \to \bbT^3$ be partially hyperbolic.
Then $f$ lifts to $\tf:\bbR^3 \to \bbR^3$ and $\tf$ is a finite distance
from a linear map $A:\bbR^3 \to \bbR^3$.
Let $\lam_i$ be the eigenvalues of $A$ ordered so that
$|\lam_1|  \le  |\lam_2|  \le  |\lam_3|$.
In this section, we consider the case $|\lam_2| = 1$
and give an outline of the proof of \cref{thm:nil}
in this setting.

Replacing $f$ by an iterate if necessary,
we assume that all of the eigenvalues are positive.
As noted at the end of \cref{sec:reebless},
the eigenvalues satisfy
$\lam_1 < 1 < \lam_3$.

If $\cF$ is a Reebless branching foliation in $\bbT^3$,
then $\tF$ is almost parallel to a linear foliation,
that is, the translates of a plane.
If $A(\tF)$ is almost parallel to $\tF$,
then this plane is invariant and is therefore
the span of two eigenspaces of $A$.
Denote the possible linear foliations as
$\tAus$, $\tAcu$, and $\tAcs$ corresponding to the
partially hyperbolic splitting of $A$.
We consider the three possibilities in turn.

\begin{prop} \label{prop:ustorus}
    If $\tFcs$ is almost parallel to $\tAus$,
    then $\Fcs$ contains a compact leaf.
\end{prop}

Since \cref{thm:nil} assumes no compact leaves, this eliminates one of the
possibilities.

\begin{idea}
    Since $\lam_2 = 1$, there is a rank-2 subgroup
    $\Gam \subset \pi_1(\bbT^3)$ corresponding to $\EsA \oplus \EuA$.
    Take any $L \in \tFcs$.
    By \cref{cor:halfspace},
    $L$ splits $\bbR^3$ into two half spaces $L^+$ and $L^-$.
    Consider the union $\bigcup_{\gam \in \Gam} \gam(L^+)$.
    Since $\tFcs$ is almost parallel to $\tAus$,
    this union is not the entire space
    $\bbR^3$ and therefore has a non-empty boundary.
    Each component of this boundary is a leaf of $\tFcs$ which is invariant
    under $\Gam$ and therefore quotients down to a 2-torus embedded in $\bbT^3$.
\end{idea}

\begin{lem} \label{lem:boxgrow}
    $\tFcu$ is not almost parallel to $\tAcs$.
\end{lem}
\begin{idea}
    Consider a coordinate system on $\bbR^3$
    such that $A$
    is given by
    \[
        A(x,y,z) = (\lam_1 x, \lam_3 y, z).
    \]
    This change of coordinates affects volumes and distances
    in $\bbR^3$ by at most a constant factor.
    Consider
    the cube
    \begin{math}
        X = [-1,1] \times [-1,1] \times [-1,1].
    \end{math}
    As $\tf$ is a bounded distance from $A$,
    one can show that there is $C>1$
    such that
    \[
        \tf^n(X) \subset
        [-C,C] \times [-(\lam_3+C)^n, (\lam_3+C)^n] \times [-C n,C n].
    \]
    for all $n>0$.
    Now suppose $J$ is an unstable curve inside $X$.
    If $\tFcu$ is almost parallel to $\tAcs$,
    then each $\tf^n(J)$ lies in an $R$-neighbourhood of a leaf of $\tAcs$.
    This means there is a sequence
    $\{y_n\}$
    such that
    \[
        \tf^n(J) \subset
        [-C,C] \times [y_n-R,y_n+R] \times [-C n,C n].
    \]
    As the volume of this box grows only linearly in $n$,
    \cref{lem:volume} gives a contradiction.
\end{idea}
This last proof is a specific case
of the idea of confining an unstable segment of $\tf$
to the $R$-neighbourhood of a leaf of a model foliation and using
\cref{lem:volume} to derive a contradiction.
This technique occurs many times in the proofs of the classification results
and we refer to it as the ``length versus volume'' argument.

Analogously,
$\tFcs$ is not almost parallel to $\tAcu$
and therefore $\tFcs$ is almost parallel to $\tAcs$.
By employing the semiconjugacy of Franks (\cref{thm:franks}),
one can show that there is a map
$\Hcs: \bbR^3 \to \tAcs$
defined by the property that
$\Hcs(p)$ is the unique leaf of $\tAcs$ for which
\[
    \sup_{n \in \bbZ} \dist(\tf^n(p), A^n \Hcs(p) ) < \infty.
\]
This map is continuous, surjective, and
commutes with deck transformations.
Further, the supremum above is bounded by a constant $S$
independent of $p$.

The most important property of $\Hcs$ is that its fibers are exactly the leaves
of $\tFcs$.
To see this, first suppose points
$p$ and $q$ lie on the same leaf $L \in \tFcs$.
By the definition of almost parallel,
$\tf^n(L)$ lies in the $R$-neighbourhood of a leaf $\cL_n \in \Acs$
and therefore
\[
    \dist(\Hcs \tf^n(p), \cL_n)  \le 
    \dist(\Hcs \tf^n(p), \tf^n(p)) + \dist(\tf^n(p), \cL_n) < S+R.
\]
The same inequality holds with $q$ in place of $p$ and so
\[
    \dist(A^n \Hcs(p), A^n \Hcs(q)) =
    \dist(\Hcs \tf^n(p), \Hcs \tf^n(q)) < 2(S+R).
\]
Then
\[
    \dist(\Hcs(p), \Hcs(q)) < 2(S+R) \lam_3^{-n}
\]
which in the limit $n \to \infty$ implies $\Hcs(p)=\Hcs(q)$.

Now suppose distinct leaves $L, L' \in \tFcs$
are both mapped by $\Hcs$ to $\cL \in \tAcs$.
A priori,
$L$ and $L'$ could be disjoint
or, in the case of a branching foliation,
could intersect without being equal.
Either way,
the region between $L$ and $L'$ has non-empty interior $U$
and $\Hcs(p) = \cL$ for every point $p \in U$.
Let $J$ be an unstable curve in $U$.
Every point in $\tf^n(J)$ is mapped to $A^n(\cL)$.
A minor adaptation of the proof of \cref{lem:boxgrow}
then gives a contradiction.

Thus the fibers of $\Hcs$ are disjoint surfaces
and form a true, non-branching foliation
tangent to $\Ecsf$.
Since the definition of $\Hcs$ did not depend on $\Fcs$,
this invariant foliation is unique.

A similar construction holds for $\Hcu:\bbR^3 \to \tAcu$
The product $H = \Hcs \times \Hcu$ yields a homeomorphism
between the space of center leaves of $\tf$ and those of $A$.
Let $P:\bbR^3 \to \bbR$ be a linear map transverse to the
center direction of $A$.
Using this $P$, the averaging method of Fuller
described in \cref{sec:leafconj}
gives the leaf conjugacy.

\section{Euclidean geometry, isotopic to Anosov} \label{sec:da} %

Let $f:\bbT^3 \to \bbT^3$ be partially hyperbolic.
Then $f$ lifts to $\tf:\bbR^3 \to \bbR^3$ and $\tf$ is a finite distance
from a linear map $A:\bbR^3 \to \bbR^3$.

Let $\lam_i$ be the eigenvalues of $A$ ordered so that
$|\lam_1|  \le  |\lam_2|  \le  |\lam_3|$.
In contrast to \cref{sec:skew}, we now discuss
the proof of \cref{thm:nil} for the case $|\lam_2|  \ne  1$.
Up to replacing $f$ by its inverse, assume $|\lam_2| > 1$

\begin{lem} \label{lem:real}
    The eigenvalues $\lam_i$ are real and distinct.
\end{lem}
\begin{idea}
    Assume $\lam_2$ and $\lam_3$ are complex.
    Then there is a unique $A$-invariant two-dimensional subspace $P \subset
    \bbR^3$ and $A \inv|_P$ is a contraction.
    By this uniqueness, $\tFcs$ is almost parallel to the translates of $P$.
    One can then consider a stable arc $J$ and
    use the contraction of $A \inv|_P$ in a length versus volume
    argument to get a contradiction.
    Thus, the $\lam_i$ are real.
    Algebraic properties of $SL(3, \bbZ)$ imply that the
    eigenvalues are distinct (see for example \cite[Definition V.3.10]{hungerford} and the remark that follows).\footnote{
    An elementary proof is as follows.  Suppose $A \in SL(3, \bbZ)$ has
    irrational eigenvalues $\lam_1 = \lam_2 < \lam_3$
    counted with multiplicity.
    Irrationality implies that $\det(x I - A)$
    is the minimal polynomial for the $\lam_i$.
    Using $\lam_1^2 \lam_3 = \det(A) = 1$ and $2 \lam_1 + \lam_3 = T$
    where $T \in \bbZ$ is the trace of $A$,
    observe that $2 \lam_1^3 - T \lam_1^2 + 1 = 0$.
    Uniqueness of the minimal polynomial implies that
    $\det(x I - A) = x^3 - \tfrac{T}{2} x^2 + \tfrac{1}{2}$
    which is impossible for a matrix with integer entries.
    Thanks goes to M.~Coons for aiding in this proof.}
\end{idea}
There are exactly three linear foliations with property
that $A(\tA)$ is almost parallel to $\tA$.
These are the foliations $\tAcs$, $\tAcu$, and $\tAus$
corresponding to the partially hyperbolic splitting.
Descending to the compact manifold,
$\Acs$, $\Acu$, and $\Aus$
are foliations of $\bbT^3$ by planes.

\begin{prop} \label{prop:planes}
    Let $\cF$ be a two-dimensional (branching) foliation on $\bbT^3$
    and $\cW$ a transverse one-dimensional foliation.
    Suppose $\cA$ is a linear foliation by planes
    such that the lifts $\tF$ and $\tA$ are almost parallel.
    Then
    \begin{enumerate}
        \item $\cF$ is a (branching) foliation by leaves diffeomorphic to planes;

        \item the lifted foliations $\tF$ and $\tW$ have global product structure;
        that is, each leaf of $\tF$ intersects each leaf of $\tW$
        in exactly one point; and

        \item there are $R > 1$ and a closed cone $\cE \subset \bbR^3$
        such that
        \[            \|p - q\| > R  \quad \Rightarrow \quad  p - q \in \cE  \]
        for any points $p$ and $q$ lying on the same leaf of $\tW$.
        Moreover, any line $E \subset \cE$ is transverse to $\tA$.
    \end{enumerate}  \end{prop}
\begin{idea}
    If $\cF$ were not a foliation by planes,
    there would be $L \in \tF$ and $\gam \in \pi_1(\bbT^3)$
    such that $\gam(L)=L$.
    One can then show that $L$ does not lie in an $R$-neighbourhood
    of any $\cL \in \tA$.
    If every leaf of a foliation is simply connected,
    the foliation is without holonomy and global product structure
    follows from results in foliation theory (see  \cite{Pot}).
    
    The third item can be proved using global product structure and the
    definition of almost parallel.
\end{idea}
Let $H:\bbR^3 \to \bbR^3$ be the lifted semiconjugacy of Franks
(\Cref{thm:franks}).
It is a finite distance from the identity and
satisfies $H \tf = A H$.

\begin{lem}
    $\tFcu$ is almost parallel to $\tAcu$.
\end{lem}
\begin{idea}
    Let $\tA$ be the linear foliation almost parallel to $\tFcu$.
    Note that
    \begin{align*}
        q \in \tW^s(p)
        & \quad \Rightarrow \quad 
        \lim_{n \to \infty} \|H \tf^n(q) - H \tf^n(p)\| = 0 \\
        & \quad \Rightarrow \quad 
        \lim_{n \to \infty} \|A^n H(q) - A^n H(p)\| = 0 \\
        & \quad \Rightarrow \quad 
        H(q) - H(p) \in \EsA
    \end{align*}
    where $\EsA$ is the eigenspace associated to $\lam_1$.
    By \cref{prop:planes},
    $\EsA$ is contained in a cone $\cE^s$ transverse to $\tA$.
    Hence, $\tA = \tAcu$
    as that is the only linear $A$-invariant foliation
    transverse to $\EsA$.
\end{idea}
\begin{lem} \label{lem:dapar}
    $\tFcs$ is almost parallel to $\tAcs$.
\end{lem}
\begin{idea}
    First, suppose $\tFcs$ is almost parallel to $\tAcu$.
    Then, the same argument as in the proof of \cref{lem:real}
    gives a contradiction

    Next, suppose $\tFcs$ is almost parallel to $\tAus$
    and let $\cE^c$ be the cone transverse to $\tAus$
    from \cref{prop:planes}.
    If $q \in \tWu(p)$,
    then $\tf^n(p) - \tf^n(q) \in \cE^c$
    for large positive $n$.
    Then,
    \[    
        H(\tf^n(p)) - H(\tf^n(q)) =
        A^n H(p) - A^n H(q) =
        A^n \bigl(H(p)-H(q)\bigr)
    \]
    is bounded away from $\EuA$ as $n \to \infty$.
    This is only possible if $H(p) - H(q) \in \EcA \oplus \EsA$.

    The previous proof showed that if $q \in \tWs(p)$,
    then $H(p) - H(q) \in \EsA$.
    Thus for any leaf $\cL \in \tAcs$,
    $H \inv(\cL)$ is saturated by both the stable and unstable
    leaves of $\tf$.
    This is a non-generic property.
    By perturbing $f$ in the space of $C^1$ diffeomorphisms,
    we may assume that $f$ is \emph{accessible}
    and consequently that $X = \bbR^3$ is
    the only non-empty subset $X \subset \bbR^3$
    which is saturated by both stable and unstable leaves
    of $\tf$. One can show that the perturbation does not modify the asymptotic directions of the foliations. 
    This gives a contradiction.

    As $\tAcu$ and $\tAus$ have been eliminated,
    the only remaining possibility is that
    $\tFcs$ is almost parallel to $\tAcs$.
\end{idea}
Using \cref{lem:dapar} and the properties of $H$,
one can show that for any $L \in \tFcs$
the image $H(L)$ is equal as a subset of $\bbR^3$
to a leaf of $\tAcs$.
Suppose distinct leaves $L, L' \in \tFcs$ satisfy
$H(L)=H(L')$.
By the global product structure given by \cref{prop:planes},
there is an unstable curve $J$ with endpoints
$p \in L$ and $q \in L'$.
For large positive $n$, $\tf^n(p) - \tf^n(q) \in \cE^u$
where $\cE^u$ is a cone transverse to $\tAcs$
given by \cref{prop:planes}.
Therefore,
$H \tf^n(p) - H \tf^n(q)$
is bounded away from $\EcsA$
as $n \to \infty$ giving a contradiction.

This shows that $H$ defines a bijection between the
leaves of $\tFcs$ and those of $\tAcs$.
A similar result holds for $\tFcu$
and the remaining steps for establishing the leaf conjugacy
are exactly as for the systems considered in \cref{sec:skew}.

\section{Nil geometry} \label{sec:nil} %

We next consider the proof of \cref{thm:nil}
for 3-manifolds with Nil geometry.
These are non-trivial circle bundles over $\bbT^2$.
As the bundle is non-trivial, there cannot be a horizontal surface.
Because of this property,
all leaves of a Reebless foliation are vertical.
In some ways, this makes these manifolds simpler to study than $\bbT^3$
even though the geometry is more difficult to visualize.

To give a precise description of the Reebless foliations,
we first look at the universal cover.
Define Heisenberg space $\Heis$ as the nilpotent Lie group
consisting of matrices of the form
\[
        \begin{pmatrix}
        1 & x & z \\ 0 & 1 & y \\ 0 & 0 & 1  \end{pmatrix}
    .
\]
In this section, we write the above element as $(x,y,z) \in \Heis$
for brevity.
Every compact 3-manifold with Nil geometry can be defined as
the quotient of $\Heis$ by a discrete subgroup.
Let $\cN$ be such a manifold.

Let $\pi:\Heis \to \bbR$ be of the form $\pi(x,y,z) = a x + b y$
for some non-zero $(a,b) \in \bbR^2$.
The level sets of $\pi$ define a foliation $\tA_\pi$ of $\Heis$ by planes
which quotients down to a foliation $\cA_\pi$ on $\cN$.

\begin{prop}
    If $\cF$ is a Reebless foliation on $\cN$, its lift $\tF$
    is almost parallel to $\tA_\pi$ for some $\pi$.
\end{prop}

This is a restatement of \cref{thm:models} in this specific setting.
Now suppose $f:\cN \to \cN$ is partially hyperbolic.
Since $\cN$ is a nilmanifold,
the results of Mal'cev \cite{Malcev} imply that there is a
Lie group automorphism $\Phi:\Heis \to \Heis$
such that $\Phi$ quotients down to a diffeomorphism $g:\cN \to \cN$
homotopic to $f$.
This automorphism is of the form
\[
    \Phi(x,y,z) = (A(x,y), z + p(x,y))
\]
where $A: \bbR^2 \to \bbR^2$ is a linear map with $\det(A) = 1$
and $p$ is a quadratic polynomial.

A length versus volume argument implies that $A$ has an eigenvalue
$\lam > 1$ and so $g$ is partially hyperbolic.
Let $\tAcs$ and $\tAcu$ be the invariant foliations of $\Phi$.
Since $\Eug \oplus \Esg$ is not integrable,
there is no $\tAus$ foliation.

For a one dimensional subspace 
$E \subset \bbR^2$, the following properties are equivalent{:}
\begin{itemize}
    \item $E$ lies at a finite distance from $A(E)$,

    \item $E$ is an eigenspace of $A$.
\end{itemize}
From this, one can show that if $\tA_\pi$ is almost parallel to
$\Phi(\tA_\pi)$, then either $\tA_\pi = \tAcs$ or $\tA_\pi = \tAcu$.
If $\tFcu$ is almost parallel to $\tAcs$,
a variant of the proof of \cref{lem:boxgrow}
gives a contradiction.
Therefore, $\tFcu$ is almost parallel to $\tAcu$
and similarly $\tFcs$ is almost parallel to $\tAcs$.
From this point on, the proof is as in \cref{sec:skew}
with minor variations.  For the averaging method of Fuller,
use $P(x,y,z) = z$.

\section{Sol geometry} \label{sec:sol} %

This section discusses the proof of \cref{thm:sol}.
Let $A:\bbT^2 \to \bbT^2$ be a hyperbolic toral automorphism
and define the suspension flow $\varphi$
on the manifold $M_A = \bbT^2 \times \bbR / (Av,s) \sim (v,s+1)$.

Suppose $f$ is a partially hyperbolic diffeomorphism
on $M_A$.
By replacing $f$ with an iterate,
we may assume $f$ is homotopic to the identity.
Then $f$ lifts to $\tf$ on the universal cover $\tM_A$ where $\tf$ is a finite
distance from the identity.

The universal cover $\tM$ is homeomorphic to $\bbR^2 \times \bbR$,
and so we will at times write points in $\tM$ in the form $(v,s)$,
but keep in mind that the geometry is not Euclidean.
The Anosov flow $\varphi$ lifts to a flow $\tvp$
given by $\tvp_t(v,s) = (v,s+t)$.
The time-one map is partially hyperbolic.
Let $\tAcs$, $\tAcu$, and $\tAus$ be the corresponding model foliations.

\begin{prop}
    Any Reebless foliation $\cF$ on $M_A$ lifts to
    a foliation $\tF$ on $\tM_A$ which is almost parallel to one of
    $\tAcs$, $\tAcu$, and $\tAus$.
\end{prop}

The following result on compact leaves is proved
in a similar way to \cref{prop:ustorus}.

\begin{prop}
    If $\cF$ lifts to a foliation almost parallel to $\tAus$,
    then $\cF$ has a compact leaf.
\end{prop}

As we are under the hypothesis that
$\tFcs$ has no torus leaves,
it follows that $\tFcs$ is almost parallel either to
$\tAcs$ or $\tAcu$.
Up to possibly replacing $A$ by $A \inv$,
we may freely assume that $\tFcs$ is almost parallel to $\tAcs$.
Our goal is then to show that $\tFcu$ is almost parallel to $\tAcu$.

\begin{lem} \label{lem:cuap}
    $\tFcu$ is almost parallel to $\tAcu$.
\end{lem}

This step was surprisingly difficult.
It was the last piece of the puzzle when the results
in \cite{HP2} were developed
and frustrated the authors for some time.
Sadly,
its proof is not easily summarized
and so for this survey we just take \cref{lem:cuap} as a given.

\begin{lem}
    For every $L \in \tFcs$,
    there is a unique $\cL \in \tAcs$
    such that $L$ and $\cL$ are at finite distance.
\end{lem}
\begin{idea}
    By the definition of almost parallel,
    at least one such leaf exists.
    Since the elements of $\tAcs$ are the weak stable leaves
    of a lifted Anosov flow, no two of them are at finite distance.
\end{idea}
There is also a converse result, swapping the roles of
$\tFcs$ and $\tAcs$.

\begin{lem}
    For every $\cL \in \tAcs$,
    there is a unique $L \in \tFcs$
    such that $L$ and $\cL$ are at finite distance.
\end{lem}
\begin{idea}
    Suppose two leaves $L,L' \in \tFcs$
    are both at finite distance from $\cL \in \tAcs$.
    Then one can take an unstable curve $J$ lying
    between $L$ and $L'$ and apply a length versus volume
    argument for $\tf^n(J)$ to derive a contradiction.
    In this case, the details are complicated and
    one must consider volume on an intermediate covering space
    homeomorphic to $\bbT^2 \times \bbR$.
\end{idea}
\begin{cor}
    $\tf(L) = L$ for all $L \in \tFcs$.
\end{cor}
\begin{idea}
    Since $\tf$ is a finite distance from the identity,
    both $L$ and $\tf(L)$ are a finite distance from the same leaf
    $\cL \in \tAcs$.
\end{idea}
Since we are trying to compare the system to a flow,
knowing that $\tf$ fixes leaves is an important step.
Paradoxically, fixing leaves implies that $\tf$ does not fix points.
This was first observed in \cite{BW}.

\begin{lem} \label{lem:noper}
    $\tf$ has no periodic points.
\end{lem}
\begin{idea}
    Suppose $p = \tf^n(p)$ is a periodic point.
    Take $q \in \tWu(p)$ such that $q$ lies on a leaf $L \in \tFcs$
    with $p \notin L$.
    Then $\tf^{-n k}(q) \to p$ as $k \to \infty$
    and $\tf^{-n k}(q) \in L$ for all $k$,
    a contradiction.
\end{idea}

Now consider a point $p \in \tM_A$.
Take leaves $\Lcs \in \tFcs$ and $\Lcu \in \tFcu$
such that $p \in \Lcs \cap \Lcu$.
By a compactness argument,
one can show that $\Lcs \cap \Lcu$ has a finite number of connected components,
each of which is a properly embedded line tangent to the center
bundle of $\tf$.
Thus, if $L$ is the component containing $p$,
there is $n  \ge  1$ such that $\tf^n(L)=L$.
As $\tf^n|_L$ is a fixed-point-free homeomorphism
of a line,
the orbit $\{\tf^{n k}(p) : k \in \bbZ \}$ is unbounded in $L$.

Using the coordinates $(v,s)$ for points in $\tM_A$,
write $\tf^{n k}(p) = (v_k, s_k)$.
Unboundedness inside $L$ implies that
$\inf s_k = -\infty$
and
$\sup s_k = +\infty$.
Expansiveness of the Anosov flow can be used
to show that there is a unique leaf $\cL \in \tAcs$
such that
\[
    \sup_{k \in \bbZ} \dist(\tf^{n k}(p), \cL) < \infty.
\]
This defines the map $\Hcs:\tM \to \tAcs$.
One can verify that the fibers of $\Hcs$ are exactly
the leaves of $\tFcs$.
Hence, $\tFcs$ is a true foliation
and is the unique invariant foliation tangent to
the center-stable direction of $\tf$.
A similar map $\Hcu:\tM \to \tAcu$ exists
and their product gives a homeomorphism
between the space of center leaves of $\tf$
and the orbits of the flow $\tvp$.
Using the averaging method of Fuller
with $P(v,s)=s$,
one can construct a leaf conjugacy.

\section{Seifert spaces and beyond} \label{sec:seifert}  %

We now discuss Reebless foliations on Seifert spaces
and give an outline of the proof of \cref{thm:phseifert}.

For simplicity, we only consider circle bundles, instead of Seifert spaces
containing exceptional fibers.
We also need only consider circle bundles over surfaces of genus $g \ge 2$,
since circle bundles over $\bbS^2$ and $\bbT^2$ have nilpotent
fundamental group.

Suppose $M$ is a circle bundle over a surface $\Sig$.
The circle bundle is given by the fibers of a map $p:M \to \Sig$.
Call an immersed surface $L$ \emph{horizontal} if $p|_L : L \to \Sig$
is a local homeomorphism.
Call $L$ \emph{vertical} if $L = p \inv(\alpha)$ where $\alpha$ is an immersed
topological circle or line in $\Sig$.

Brittenham proved the following result \cite{Brittenham}.

\begin{thm} \label{thm:brittenham}
    Let $M$ be a circle bundle over a surface $\Sig$
    and
    let $\cF$ be a Reebless foliation on $M$.
    Then, the bundle map $p: M \to \Sig$ may be chosen so that
    every leaf of $\cF$ is either horizontal or vertical.
\end{thm}

The regularity of the foliation is crucial.
\Cref{thm:brittenham}
builds on unpublished work of Thurston, who showed that if the Reebless
foliation $\cF$ is $C^2$, then every leaf is horizontal.
See \cite{Levitt} for details.
In $C^1$ regularity or lower, it is possible to have a mixture of horizontal
and vertical leaves.
See \cite[Chapter 4]{Calegari} for examples, as well as an introduction
to subject and a simplified outline of the proof of \cref{thm:brittenham}.

In the partially hyperbolic setting,
the branching foliation $\Fcs$ or its approximation $\Fcsep$ is $C^0$
with $C^1$ leaves.
However, we may use the dynamics to show that a mix of horizontal and
vertical leaves is impossible.

\begin{prop} \label{prop:novert}
    Let $f$ be a transitive partially hyperbolic diffeomorphism
    defined on a circle bundle $M$ over a surface $\Sig$ of genus $g  \ge  2$.
    Then every leaf of $\Fcs$ is horizontal.
\end{prop}

\begin{idea}
    For simplicity, assume $f$ is dynamically coherent,
    so that $\Fcs$ is a true foliation.
    Applying \cref{thm:brittenham},
    we may consider each leaf of $\Fcs$ as either horizontal or vertical.
    Let $\Lam \subset M$ be the union of all vertical leaves in $\Fcs$.
    Note that $\Lam  \ne  M$, since if every leaf of a foliation on $M$ were
    vertical, it would imply that there is a one-dimensional foliation
    on the higher genus surface $\Sig$.
    To show there are no vertical leaves we assume $\Lam \ne \varnothing$
    and derive a contradiction.

    Let $\delta > 0$ be small and let $W^u_\delta(x)$ denote the unstable
    segment of length $\delta$ centered at $x \in M$.
    Since $f$ maps vertical leaves to vertical leaves, $f(\Lam) = \Lam$.
    One can also show that $\Lam$ is a closed subset of $M$.
    If $\varnothing  \ne  \Lam  \ne  M$, then there is $\delta > 0$ such that
    $U := \bigcup_{x \in \Lam} \Ws_\delta(x)$
    satisfies
    $\varnothing  \ne  U  \ne  M$.
    Then, $f^{-n}(\overline U) \subset U$ for all large $n$ and this contradicts
    transitivity.
      \end{idea}

\begin{proof}
    [Idea of the proofs of \cref{thm:phseifert} and \cref{cor:noproduct}.]
    Since all of the leaves of $\Fcs$ are horizontal,
    they are transverse to the fibers of the Seifert fibering.
    In particular, no unit vector in $\Ecs$ is tangent to the fibering.
    Assuming that $\Ec$ is orientable, this gives us a vector field $v$
    which is nowhere tangent to the fibering.
    Assume that the bundle map $p:M \to \Sig$ has derivative
    $D p:TM \to T\Sig$
    and so the vector field $v:M \to TM$ composed with $D p$
    gives a non-zero map $\rho:M \to T\Sig$.

    Note that at this point we can prove \cref{cor:noproduct} directly.
    Suppose $M$ is the direct product $\Sig \times \bbS^1$.
    Then, there is a section $\Sig \hookrightarrow M$.
    For instance, map $\Sig$ to $\Sig \times \{1\}$.
    Composing the map $\rho$ with this section yields
    a non-vanishing vector field on $\Sig$.
    Since $\Sig$ is a higher-genus surface, this is impossible.

    We now go back to the general case of \cref{thm:phseifert}.
    Our goal is to show that $M$ supports an Anosov
    flow.
    To do this, one uses $x \mapsto \rho(x)/\|\rho(x)\|$
    to define a map
    from $M$ to the unit tangent bundle $T^1\Sig$
    taking fibers to fibers.
    The properties of circle bundles imply that there is a finite covering map
    from $M$ to $T^1\Sig$.
    Then, the geodesic flow on $T^1\Sig$ lifts to an Anosov flow on $M$.
\end{proof}

At the time of writing, it is not clear if the assumption of
transitivity is necessary in \cref{thm:phseifert}.
It seems possible that a careful analysis of the possible dynamics
on vertical $cs$ and $cu$-leaves could be used to establish
\cref{thm:phseifert} for all partially hyperbolic systems.

After Seifert spaces, the next ``natural'' family of manifolds to consider
are ones with hyperbolic geometry (or atoroidal).

Goodman demonstrated that some hyperbolic 3-manifolds support Anosov flows \cite{goodman}, see also \cite{Fenley}.

There is an established history of research for Reebless foliations
on hyperbolic 3-manifolds \cite[Chapters 7 and 8]{Calegari}.
However, the study of partially hyperbolic systems on
these manifolds has only barely begun \cite[Section 3]{Parwani}.
For instance, many hyperbolic 3-manifolds, such as the Weeks manifold, do not
support Reebless foliations and so do not support partially hyperbolic
diffeomorphisms \cite{CalegariDrunfield,RSS}.

We refer the reader to \cite{barbotfenley-toroidal} and references therein for results concerning Anosov and pseudo-Anosov flows in 3-manifolds which are not atoroidal.

\section{Other classification results} \label{sec:others} %

So far, we have only considered restrictions on the manifold $M$.
We now consider assumptions on the dynamics itself.
We first discuss a very strong condition for the invariant foliations
called quasi-isometry.

Consider a foliation $\cF$ on a manifold $M$.
Lift $\cF$ to $\tF$ on the universal cover.
For two points $x,y$ on the same leaf $L \in \tF$
one can define a distance $d_{\tF}(x,y)$
as the length of the shortest path inside the leaf $L$.
Ignoring $L$, one can also define $d_{\tM}(x,y)$
as distance on the ambient manifold $\tM$.
The foliation $\cF$ is \emph{quasi-isometric}
if these two notions of distance are roughly the same;
to be precise, if there is a uniform constant $Q > 1$
such that
\[
    d_{\tF}(x,y) < Q d_{\tM}(x,y) + Q
\]
for any $x,y \in \tM$ lying on the same leaf of $\tF$.

The use of the word ``quasi-isometric''
here is slightly different than the much more common notion
of two metric spaces being quasi-isometric to each other. 

The concept of a quasi-isometric foliation was first applied
to partially hyperbolic systems by Brin,
who proved the following \cite{Brin}.

\begin{thm} \label{thm:qidc}
    Suppose $f:M \to M$ is absolutely partially hyperbolic
    and that each of $\Wu$ and $\Ws$ is a quasi-isometric foliation.
    Then, $\Ecs$, $\Ecu$, and $\Ec$ are uniquely integrable.
\end{thm}
\begin{idea}
    If $\Ecs$ is not uniquely integrable,
    then there is a path $\alpha:[0,1] \to \tM$ piecewise
    tangent to $\Ecs$ with endpoints
    $\alpha(0)  \ne  \alpha(1)$
    which lie on the same unstable leaf.
    Since $f$ is \emph{absolutely} partially hyperbolic,
    there are constants $\mu > \gam > 1$ such that
    \[
        \|Df v^{cs} \| < \gam < \mu < \| Df v^u \|
    \]
    for all unit vectors $v^{cs} \in \Ecs$ and $v^u \in \Eu$.
    Quasi-isometry and the expansion along $\Eu$ imply that
    \[
        d_{\tM}(\tf^n \alpha(0), \tf^n \alpha(1)) >
        {\mu^n}\,{Q \inv}\, d_u(\alpha(0), \alpha(1)) \, - \, 1
    \]
    for all $n$.
    Since $\alpha$ is piecewise tangent to $\Ecs$,
    \[
        d_{\tM}(\tf^n \alpha(0), \tf^n \alpha(1)) <
        \length(\tf^n \circ \alpha) <
        \gam^n \length(\alpha).
    \]
    For large $n$, these two estimates give a contradiction.
    \begin{figure}
        [t]
        \begin{center}
            \includegraphics{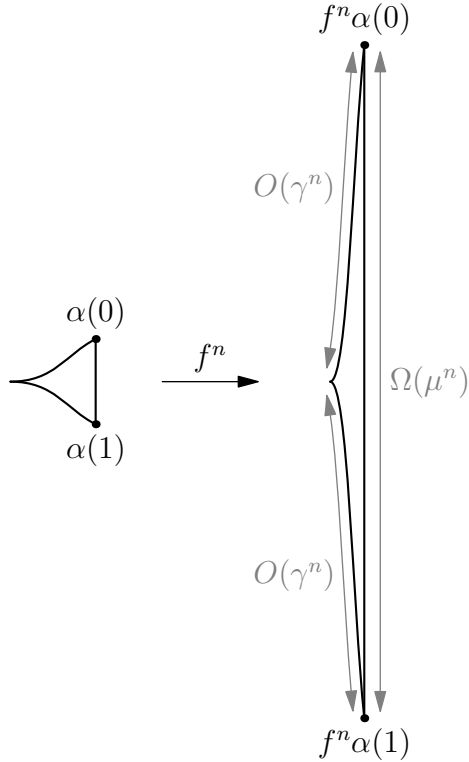}  \end{center}
        \caption{
        A depiction of the proof of \cref{thm:qidc}.
        The points $\alpha(0)$ and $\alpha(1)$ lie on a piecewise $C^1$ curve
        tangent to $\Ecs$, shown here as two curved segments joined in a
        cusp.
        They also lie on the same unstable curve, shown here as a vertical line.
        By quasi-isometry of $\Wu$,
        the distance between the endpoints grows at least on the order of
        $\mu^n$, but the length of the $\Ecs$ curve grows only on the order of
        $\gam^n$, a contradiction.
        }
        \label{fig:exp}
    \end{figure}
\end{idea}
We have given an outline of this proof not only because it is short
and elegant,
but because it highlights why it is necessary to require
absolute instead of pointwise partial hyperbolicity
when using the assumption of quasi-isometry.
In some sense, the proof is a global version of the local arguments used in
establishing \cref{Teo-StrongUnstableManifold}.

Brin's motivation in proving \cref{thm:qidc}
was to answer the question of dynamical coherence on the 3-torus.
Working with Burago and Ivanov \cite{BBI2},
he achieved this goal in the end.

\begin{thm} \label{thm:torusdc}
    Let $f:\bbT^3 \to \bbT^3$ be absolutely partially hyperbolic.
    Then, each of $\Wu$ and $\Ws$ is a quasi-isometric foliation.
    As a consequence, $f$ is dynamically coherent.
\end{thm}
This theorem was later extended by Parwani to all nilmanifolds in dimension 3
\cite{Parwani}.
It was also a major motivation for the classification up to leaf conjugacy
of partially hyperbolic systems on $\bbT^3$.
In fact,
a version of the classification holds for higher dimensional tori
when quasi-isometry is assumed.

\begin{thm} [\cite{Hammerlindl}] \label{thm:hightorus}
    Suppose that $f:\bbT^d \to \bbT^d$ is absolutely partially hyperbolic
    with $\dim \Ec = 1$ and that each of $\Ws$ and $\Wu$
    is a quasi-isometric foliation.
    Then $f$ is leaf conjugate to a linear automorphism of $\bbT^d$.
\end{thm}
Many dynamical properties can be established when assuming
quasi-isometry of the foliations
\cite{Hammerlindl, ham-dynqi, ham-pgps}.
Quasi-isometry is such a strong assumption, however,
that it is unlikely to hold for the vast majority
of partially hyperbolic systems of interest.
Brin briefly discusses this in \cite{Brin}.
Several conditions which prevent $\Ws$ and $\Wu$ from being
quasi-isometric are given in \cite{ham-expfoln}. 

The quasi-isometric condition and absolute domination in Theorem \ref{thm:hightorus} can be removed if one assumes that $f$ is isotopic to a linear Anosov automorphism by a path of partially hyperbolic diffeomorphisms (see \cite{FPS} for precise and more general statements, in particular, in dimension 3, stronger results hold \cite{Pot}).  In \cite{potrietrapping} conditions on the center-stable foliation are given that make the diffeomorphism resemble a linear Anosov automorphism of a torus or a nilmanifold.

We suspect that the 3-torus and other nilmanifolds
are the only 3-di\-men\-sion\-al manifolds where
the quasi-isometry property holds for partially hyperbolic systems.
If one also considers quasi-isometry for the center foliation,
then the following holds.

\begin{thm}
    [\cite{ham-pgps}]
    For an absolutely partially hyperbolic system on a 3-manifold,
    the stable, center, and unstable foliations exist and are quasi-isometric
    if and only if
    $M$ is finitely covered by the 3-torus.
\end{thm}

Bonatti and Wilkinson \cite{BW} investigated the properties of transitive partially
hyperbolic systems in dimension 3.
Their original plan was to construct new examples of such systems
by gluing together pieces of already known examples.
For instance,
they hoped to construct a system which resembled
a skew product in one part of the manifold
and the time-1 map of a flow elsewhere.
They soon realized that such surgery techniques were not possible
in this setting and proved the following theorems.

\begin{thm} \label{thm:transkew}
    Let $f$ be a partially hyperbolic transitive diffeomorphism on
    a 3-manifold $M$.
    Assume that $\gam=f(\gam)$ is a circle tangent to $\Ec$
    and that $\Ws(\gam) \cap \Wu(\gam)$
    contains at least one other connected component which is a circle.
    Then, $f$ is finitely covered by a skew product.
\end{thm}

\begin{thm} \label{thm:transflow}
    Let $f$ be a dynamically coherent, transitive diffeomorphism
    on a 3-manifold.
    Assume $\gam=f(\gam)$ is a circle tangent to $\Ec$
    such that $f(L)=L$ for every center leaf $L \subset \Ws(\gam)$.
    Then $f$ is topologically conjugate
    to the time-1 map of a continuous expansive flow.
\end{thm}
Further consequences also hold.
See \cite{BW} for the precise
results.  It is further conjectured in \cite{BW} that the continuous expansive flow obtained in the above theorem is orbit equivalent to an Anosov flow. In our view, showing leaf conjugacy to such an expansive flow would be enough for the purposes of the classification of partially hyperbolic diffeomorphisms.

The proofs of the above two theorems rely on the following property:
if $f$ is transitive and $\gam$ is a periodic compact center leaf,
then each of $\Wu(\gam)$ and $\Ws(\gam)$ is dense in $M$.
The intersections of $\Wu(\gam)$ and $\Ws(\gam)$
may then be analyzed to establish global properties of the dynamics.

Say that a partially hyperbolic diffeomorphism has a
compact center foliation
if there is an invariant foliation tangent to $\Ec$
consisting entirely of compact leaves.
One step in proving \cref{thm:transkew} is the following result,
which does not rely on transitivity.

\begin{thm} \label{thm:skew3}
    If $f$ is a partially hyperbolic system
    with compact center foliation on a 3-manifold $M$,
    then $f$ is finitely covered by a skew product.
\end{thm}
This result was generalized independently and at nearly the same time
by D.~Bohnet \cite{bohnet2015quotient},
P.~Carrasco \cite{carrasco2011compact}, and
A.~Gogolev \cite{gogolev2011compact}.

\begin{thm}
    If $f$ is a partially hyperbolic system
    with compact center foliation and $\dim \Eu = \dim \Es = 1$,
    then $f$ is finitely covered by a skew product.
    
\end{thm}
Gogolev also proved the following, which applies to systems in dimensions 4
and 5.

\begin{thm} \label{thm:skew212}
    If $f$ is a partially hyperbolic system
    with compact center foliation, $\dim \Ec = 1$,
    $\dim \Eu  \le  2$, and $\dim \Es  \le  2$,
    then $f$ is finitely covered by a skew product.
\end{thm}

An open question, attributed in \cite{HHU-Survey} to C.~C.~Pugh,
asks if this always holds, regardless of the dimensions
of the subbundles.

\begin{quest}
    Is every partially hyperbolic system with compact center foliation
    finitely covered by a skew product?
\end{quest}

We say that a compact foliation is
\emph{uniformly compact}
if there is a uniform bound on the volume of every leaf.
This extra condition is not redundant.
For instance, Sullivan constructed a foliation of a 5-manifold
by circles with no uniform bound on the
length \cite{sullivan1976counterexample}.
In the partially hyperbolic setting, however,
it is an open question if any compact center foliation
is also uniformly compact.
Using uniform compactness as a hypothesis,
D.~Bohnet showed the following.

\begin{thm} \label{thm:skewu1}
    If $f$ is a partially hyperbolic system
    with uniformly
    compact center foliation and $\dim \Eu = 1$,
    then $f$ is finitely covered by a skew product.
\end{thm}
Note that there is no restriction here on the dimensions of $\Es$ or $\Ec$.
See \cite{carrasco2011compact, gogolev2011compact, bohnet2015quotient}
for further results on compact center foliations.

\section{New examples}\label{ss-examples} %

Examples of partially hyperbolic diffeomorphisms (in any dimension) appear in four basic forms.

\begin{enumerate}
\item Algebraic and geometric constructions,
    including linear automorphisms of tori and nilmanifolds
    and geodesic and frame flows;
\item Skew-products, 
\item Surgery on existing examples, 
\item Deformation and composition preserving cone-fields.
\end{enumerate}

None of the four mechanisms is completely understood. Algebraic examples are not even completely understood for the construction of Anosov diffeomorphisms. Skew products where the base is more hyperbolic than the fiber can be said to be well understood, but when they work the other way around (the expansion and contraction are seen in the fibers), this is just starting to be studied and several exciting examples are starting to appear (see \cite{FarrellGogolev,GORH}). 

Surgery constructions are only partly understood for Anosov flows (see \cite{FW,HandelThurston,goodman,Fenley,BBY}).
For partially hyperbolic diffeomorphisms, we only know about obstructions to applying surgery techniques \cite{BW}. 

This section gives an overview of new examples of partially
hyperbolic diffeomorphisms in 3-manifolds appearing in \cite{BPP,BGP,BGHP}. They exploit the fourth mechanism, which in a certain sense generalizes the fact that partially hyperbolic diffeomorphisms form an open set in the space of diffeomorphisms. It depends on the possibility of having a good control on the position of the cone-fields around the invariant bundles. This is the mechanism we will present in this section. The main purpose is to construct new examples starting with time one maps of Anosov flows but such that they are not isotopic to the identity.

\subsection{Mechanism}

The key mechanism for constructing new examples is the following simple proposition: 

\begin{prop}\label{prop-mecanism}
 Let $f: M \to M$ be a partially hyperbolic diffeomorphism with splitting $TM = \Es \oplus \Ec \oplus \Eu$. Let $h: M \to M$ be a smooth diffeomorphism such that for every $x\in M$ one has: 
\begin{equation}\label{eq:transv}
 Dh(\Es(x)) \trans \Ecu(h(x)) \ \text{and } \ Dh(\Eu(x)) \trans \Ecs(h(x)) 
 \end{equation}
 \noindent then, there exists $n>0$ such that $f^n \circ h$ is partially hyperbolic.  
\end{prop}

This proof essentially follows from the classical cone-field criteria. Let us briefly give a sketch of the proof: 

\begin{idea}Let us first show that there exists $n$ such that $f^n \circ h$ preserves an unstable cone-field. 

To do this, consider first a given cone-field $\cE^u$ for $f$, that is, $Df(\overline{\cE^u(x))}) \en \cE^u(f(x))$ and vectors in $\cE^u$ are expanded by $Df$. Notice that by considering the cone-fields $Df^k (\cE^u)$ one can assume that $\cE^u$ is as narrow as one wishes. 

By compactness and the fact that $Dh(\Eu(x))$ is transverse to $\Ecs(h(x)) $ one can assume that $\cE^u$ is such that $Dh(\overline{\cE^u(x)})$ is also transverse to $\Ecs(h(x)) $. In particular, there exists $n>0$ such that $Df^n(Dh(\overline{\cE^u})) \en \cE^u$. If $n$ is large enough one can also assume that $Df^n \circ Dh$ expands vectors in $\cE^u$. 

The existence of a stable cone-field is almost symmetric (but not exactly since composition with $h$ is on the right). Indeed, by considering the cone-field $Dh^{-1}(\cE^s)$ where $\cE^{s}$ is a stable cone-field for $f$ one obtains invariance under $D(f^n \circ h)^{-1}$. 

Both cone-fields imply partial hyperbolicity by classical results in the theory (see for example \cite[Section 2.2]{CrovisierPotrie}). 
\end{idea}

\begin{obs} The same statement holds for absolute partial hyperbolicity.
\end{obs}

If $M$ is a 3-manifold whose fundamental group has exponential growth, then it cannot admit an Anosov diffeomorphism or a skew-product.
If $M$ admits an Anosov flow, then the fundamental group must have exponential
growth.  Therefore, to find new examples of partially hyperbolic
diffeomorphisms it is enough to start with an Anosov flow $\varphi_t$ of a
3-manifold $M$ and find a diffeomorphism $h$
which verifies the conditions of \cref{prop-mecanism}
and such that no iterate of $h$ is isotopic to the identity.

\subsection{Non-transitive examples}

Consider an Anosov flow $\varphi_t : M \to M$ transverse to a torus $T \subset M$. One example of such a flow is the suspension of an Anosov diffeomorphism of the torus,
but wilder examples exist (see \cite{BBY} and references therein) and it will be those that will interest us. The reason is the following: We plan to compose the time $N$-map of $\varphi_t$ with a Dehn twist along a neighborhood of $T$ so that conditions \eqref{eq:transv} are satisfied but also so that the Dehn twist is not isotopic to the identity; when considering the suspension of an Anosov diffeomorphism of the torus, this is not the case as was proved in \cite{HP2}.

By Dehn twist we mean that we consider an embedded submanifold $D$ of $M$ diffeomorphic to $[0,1] \times \TT^2$ and an integer vector $v \in \ZZ^2$ and consider the diffeomorphism of $M$ which is the identity outside $D$ and in $D$ (with the coordinates $[0,1] \times \TT^2$) is given by $(t,x) \mapsto (t, x + \rho(t) v \text{ mod }\ZZ^2) $ where $\rho:[0,1] \to [0,1]$ is a smooth function which equals $0$ in a neighborhood of $0$ and $1$ in a neighborhood of $1$.

The easiest example that will work for us was constructed by Franks and Williams (\cite{FW}) and is an Anosov flow $\varphi_t$ which is transverse to a torus $T$ and has the property that $\varphi_t(T) \cap T= \varnothing$ for every $t >0$. This implies in particular that the flow $\varphi_t$ is not transitive. 

The Dehn twist is constructed so that it is supported on a set of the form $\mathcal{U}_N = \bigcup_{0\leq t \leq N} \varphi_t(T)$ for very large $N$ so that the effect of the modification can be made to be almost negligible at the level of the derivative, this allows the twist to satisfy conditions \eqref{eq:transv}.
The details can be found in \cite{BPP} but let us say a few more words on how this is done:
\begin{itemize}
    \item Parametrize $\mathcal{U}_N$ by a diffeomorphism
        $H_N: T \times [0,1] \to \mathcal{U}_N$
        defined  by $H_N (x,\frac{t}{N}) =  \varphi_{t}(x).$
    Since the center-stable foliation is invariant under the flow,
    preimages of its leaves under $H_N$ are sets of the form $\gam \times [0,1]$
    where $\gam$ is a curve in $T$.
    Therefore, the center-stable foliation on $M$ defines a one-dimensional foliation
    $\mathcal{T}^{cs}$ on $T$.
    Similarly, a transverse foliation $\mathcal{T}^{cu}$ exists for the center-unstable
    direction.
\item If the example is done carefully, it is possible to construct a path of diffeomorphisms $\{\eta_s\}_{s \in [0,1]}$ of $T$ which is not homotopically trivial and such that
    $\eta_s(\mathcal{T}^{cs})$ is transverse to
    $\mathcal{T}^{cu}$
    and
    $\eta_s(\mathcal{T}^{cu})$ is transverse to
    $\mathcal{T}^{cs}$
    for all $s$. One requires the path to start and end in the identity and such that for every $x\in T$ one has that the curve $s \mapsto \eta_s(x)$ is homotopically non-trivial. 
\item By choosing a very large $N$, one has that the strong bundles are very close to being tangent to the tangent spaces of the one dimensional foliations $\mathcal{T}^{cs} \times \{t\}$,$\mathcal{T}^{cu} \times \{t\}$ in the $T\times [0,1]$ coordinates.  Therefore, in these coordinates, one can construct a Dehn twist $(x,s)\mapsto (\eta_s(x),s)$ which, when sent back to $\mathcal{U}_N$ via $H_N$ will verify the conditions \eqref{eq:transv}. 
\end{itemize}

By general 3-manifold topology reasons, one can show as desired that this Dehn twist will not be homotopic
to the identity.
    In \cite{BPP} a specific example is treated where one can show this via elementary methods (showing that the action in homology is non-trivial). Recently, in \cite{BZ} it is shown that this construction can be made to any non-transitive Anosov flow. 

\subsection{Transitive examples}

The construction in the previous subsection relied on the fact
    that $\varphi_t(T)$ is disjoint from $T$ for all $t > 0$.
    This allowed us to spread out the perturbation and obtain the transversality conditions. If one wishes to make a transitive example starting from a transitive Anosov flow, this condition will never be satisfied.

One of the main constructions in \cite{BGP} consists in making a quantitative version of the previous argument, applied to certain transitive Anosov flows transverse to tori (see \cite{BBY} for plenty of such examples) that after huge finite lifts will be treated in a similar way as the non-transitive case. This implies solving several technical difficulties, but at the end it is possible to do so without so many differences as in the previous case.  

To guarantee that after perturbation the new diffeomorphism is transitive, one could argue directly with some attention to be able to use the transitivity of the Anosov flow and the control on the perturbation, but there is a way to do this with an easy trick that moreover allows us to obtain even better results (construct stably ergodic examples). The key point is that if the Anosov flow is volume preserving and the Dehn twist with which we compose preserves the same volume form, the new example will be also conservative, and therefore easier to show it can be perturbed to be robustly transitive (and stably ergodic). To do this, in \cite[Section 3]{BGP} the example is constructed starting with the Anosov flow constructed in \cite{BL} and as explained before, considering a big finite cover.

\begin{obs} 
In a work in progress it is shown that it is possible to change the construction in order to avoid the huge finite lift but we will not explain this here. 
\end{obs}

In \cite{BGP} there are examples of different nature also, constructed via Dehn twists along tori, but which are not transverse to the Anosov flow. The transversality conditions can be guaranteed by ``geometric'' reasons. We will explain these examples in the following subsection.

\subsection{New examples on unit tangent bundles} 

Consider a hyperbolic surface $S$ with Riemannian metric $g$ and its geodesic flow $\varphi_t: T^1 S \to T^1 S$. If the curvature of $g$ is everywhere negative, then $\varphi_t$ is an Anosov flow \cite{Anosov}. 

Let $f: S \to S$ be any smooth diffeomorphism. The projectivization $Pf : T^1 S \to T^1 S$ defined as $Pf( v) = \frac{Df v}{\|Df v\|}$ conjugates the geodesic flow of the metric $g$ with the one obtained by the pullback $f^\ast g$. That is, if $\psi_t$ denotes the geodesic flow of $f^\ast g$ one has that: 

\[ Pf \circ \varphi_t  = \psi_t \circ Pf \]

In particular, since the conjugacy is smooth, it follows that the derivative of $Pf$ sends the invariant bundles of $\varphi_t$ onto the ones of $\psi_t$. So the key remark is:

\begin{prop}\label{prop-closemetric}
It is possible to construct $f: S \to S$ not isotopic to the identity in $S$ so that $f_\ast g$ and $g$ are very close to each other (in the $C^\infty$-topology) and with $g$  having constant negative curvature. 
\end{prop}

This is done in \cite[Section 2]{BGP} via considering a hyperbolic metric of constant curvature $-1$ with a very short closed geodesic
and choosing $f$ to be a Dehn twist along this curve, since the geodesic is short, it follows that, in the universal cover, the lift of the diffeomorphism $f$ is very close to the identity. 

It follows that the diffeomorphism $Pf \circ \varphi_N$ for large $N$ will be partially hyperbolic thanks to Proposition \ref{prop-mecanism} and not isotopic to the identity because $Pf$ is not. Since $Pf$ does not preserve the length along fibers,
it will not preserve the volume form on $T^1S$.
However, it is possible to find a diffeomorphism $C^1$-close
to $Pf$ which does preserve this volume form
and for which \cref{prop-mecanism} still applies.

In principle, this idea and the examples presented in the previous section suggest that one might be able only to make examples isotopic to Dehn twists. These have the particularity that their action in homology only has eigenvalues of modulus $1$. The proof of Proposition \ref{prop-closemetric} as explained above, requires a small geodesic along which one performs a Dehn twist. If one would like to construct $f$ which belongs to a mapping class having eigenvalues larger than $1$, then one would need to compose with crossing Dehn twists, but if one curve is short, the transverse one should be long. 

In \cite{BGHP} this problem is solved by considering maps and geodesic flows in non-connected surfaces and choosing the initial diffeomorphisms not to be the geodesic flows, but to permute the connected components and apply the geodesic flow so that the diffeomorphism visits every connected component. In this way, it is possible to modify the metrics little by little in order to apply Dehn twists in different connected components so that a sufficiently large iterate will observe the multiple Dehn twists and the return map to a connected component will be partially hyperbolic. 
In particular, it is shown: 

\begin{thm}[\cite{BGHP}]\label{teo-unittangentbundle} 
Given a surface $S$ of genus $g\geq 2$ and a diffeomorphism $f: S \to S$, there exists a (stably ergodic) partially hyperbolic diffeomorphism $F: T^1 S \to T^1 S$ isotopic to $Pf$. 
\end{thm}

\subsection{Some considerations}

The mechanism for constructing new examples needs to start with a partially hyperbolic diffeomorphism on $M$. In particular, this can only give examples on new isotopy classes, but not on different manifolds.

\begin{quest}
 Is there a manifold $M$ whose fundamental group has exponential growth which admits a partially hyperbolic diffeomorphism and does not admit an Anosov flow? 
\end{quest}

In \cite{HP2,HPS-New} one can find results in the direction of a negative answer. However, the fact that we do not understand completely Anosov flows in 3-manifolds makes this question difficult in principle.

Also, these new examples provide different phenomena in different isotopy classes. So, it is natural to ask:

\begin{quest}
 Let $f: M \to M$ be a (dynamically coherent) partially hyperbolic diffeomorphism isotopic to the identity. Is $f$ leaf conjugate to an Anosov flow?. 
\end{quest}

In this direction recall \cref{thm:transflow} stating that if $f$ is transitive and every center-leaf is fixed by $f$ and there is a circle leaf then this must be essentially the case (see the discussion after the statement of \cref{thm:transflow}). See \cite{BGP} for more questions arising naturally from these new examples. 

\medskip

\acknowledgement
The authors would like to thank
Christian Bonatti,
Sylvain Crovisier,
Andrey Gogolev,
Enrique Pujals,
Mart\'in Sambarino,
and
Amie Wilkinson
for helpful conversations.

\end{document}